\documentclass{article}
\usepackage{hyperref}
\usepackage{microtype}
\usepackage{graphicx}
\usepackage{subfigure}
\usepackage{booktabs} 

\usepackage[utf8]{inputenc}
\usepackage{tabularx}
\usepackage[table]{xcolor}
\usepackage{colortbl}

\usepackage[accepted]{icml2023}

\usepackage{array}
\usepackage{longtable}
\usepackage{caption}
\usepackage{graphicx}
\usepackage{pgfplots}
\usepackage{cancel}

\usepackage{amsmath}
\usepackage{amssymb}
\usepackage{mathtools}
\usepackage{amsthm}

\usepackage[capitalize,noabbrev]{cleveref}
\usepackage{hyphenat}
\usepackage{enumitem}

\setlist[enumerate]{leftmargin=*,noitemsep,topsep=0pt}
\setlist[itemize]{leftmargin=*,noitemsep,topsep=0pt}

\usepackage{tikz}
\pgfplotsset{compat=1.17}

\icmltitlerunning{Testing Quantum and Simulated Annealers on the Drone Delivery Packing Problem}

\begin{document}

\twocolumn[
\icmltitle{Testing Quantum and Simulated Annealers on the Drone Delivery Packing Problem}

\begin{icmlauthorlist}
\icmlauthor{Sara Tarquini}{gssi}
\icmlauthor{Daniele Dragoni}{qcrl,di}
\icmlauthor{Matteo Vandelli}{qcrl}
\icmlauthor{Francesco Tudisco}{gssi,sch}
\end{icmlauthorlist}

\icmlaffiliation{gssi}{Gran Sasso Science Institute, Viale Francesco Crispi 7, L'Aquila, Italy}
\icmlaffiliation{qcrl}{Quantum Computing Research Laboratory, Leonardo S.p.A., Via R. Pieragostini 80, Genova, Italy}
\icmlaffiliation{di}{Digital Infrastructures, Leonardo S.p.A., Via R. Pieragostini 80, Genova, Italy}
\icmlaffiliation{sch}{School of Mathematics and Maxwell Institute, University of Edinburgh, Peter Guthrie Tait Road, EH9 3FD, Edinburgh, UK}

\icmlcorrespondingauthor{Sara Tarquini}{sara.tarquini@gssi.it}
\icmlcorrespondingauthor{Daniele Dragoni}{daniele.dragoni@leonardo.com}
\icmlcorrespondingauthor{Matteo Vandelli}{matteo.vandelli.ext@leonardo.com}
\icmlcorrespondingauthor{Francesco Tudisco}{f.tudisco@ed.ac.uk}

\vskip 0.3in
]

\printAffiliationsAndNotice{} 

\begin{abstract}
Using drones to perform human-related tasks can play a key role in various fields, such as defense, disaster response, agriculture, healthcare, and many others. The drone delivery packing problem (DDPP) arises in the context of logistics in response to an increasing demand in the delivery process along with the necessity of lowering human intervention. The DDPP is usually formulated as a combinatorial optimization problem, aiming to minimize drone usage with specific battery constraints while ensuring timely consistent deliveries with fixed locations and energy budget. In this work, we propose two alternative formulations of the DDPP as a quadratic unconstrained binary optimization (QUBO) problem, in order to test the performance of classical and quantum annealing (QA) approaches. We perform extensive experiments showing the advantages as well as the limitations of quantum annealers for this optimization problem, as compared to simulated annealing (SA) and classical state-of-the-art commercial tools for global optimization. 
\end{abstract}

\section{Introduction}
Combinatorial optimization problems are widely used to model complex decision-making processes involving a large number of binary choices. Due to the combinatorial nature of the problem, the cost of computing global optima can scale exponentially in the number of variables becoming quickly very prohibitive. While conventional approaches can be effective in specific cases, it is also evident the need for more general techniques, adaptable to a bigger range of applications. In this sense, a potentially advantageous approach is represented by the use of quantum optimization techniques. Motivated by the recent progress in quantum annealing hardware for solving quadratic unconstrained binary optimization (QUBO) problems, in this work we aim to test the performance of classical and quantum annealing algorithms for the DDPP. 

Solving the DDPP using quantum annealers involves transforming it into a QUBO problem. 
When doing so, however, adhering to the hardware constraints of the quantum machine, especially qubit count and connectivity, poses a substantial challenge. 
In order to formulate DDPP as a QUBO problem, we first transform the constrained quadratic optimization problem into an integer linear programming problem (ILPP), through the addition of slack variables. Subsequently, the conversion to a QUBO is done by introducing quadratic penalty terms that equal zero for feasible solutions, and a positive quadratic penalty for infeasible solutions. 

However, if implemented naively, this relatively standard procedure leads to a prohibitive number of variables and, in turn, a prohibitive qubit count, after the embedding process onto the quantum hardware. Thus, we propose a relaxed QUBO alternative formulation that equivalently solves DDPP but with significantly fewer slack variables (and thus fewer qubits). In order to solve the resulting QUBO formulations, we employ simulated and quantum annealing samplers, using the ultimate QPU model by D-Wave: Advantage system 4.1.

Advantage QPUs can hold inputs that are almost three times larger, on average, than those that could run on previous-generation D-Wave 2000Q QPUs. The last model features 5627 qubits and 15 couplers per qubit, for a total of 35,000 couplers \cite{adv}.

We studied the performance of the two approaches and the quality of the corresponding solutions. 
Overall, our main contributions are as follows:
\begin{itemize}
    \item We propose two QUBO formulations of the DDPP, one following an established all-purpose approach and another one tailored to the specific DDPP. We show that the latter QUBO formulation requires significantly fewer slack variables than the former.
    \item Using the proposed QUBO formulation, we perform extensive experimental evaluation providing statistics on the performance, concerning problem size and chain strength, and in terms of time-to-solution, memory cost, and solution quality of D-Wave's ultimate QPU as compared to classical (simulated) annealing and global-optimization deterministic baseline approaches.  
\end{itemize}

Motivated by the rapidly progressing technological advances realizing bigger and more connected quantum devices, the overarching goal of this work is to identify, using DDPP as a prototype of a complex combinatorial optimization problem, the advantages and limitations of modern quantum annealing hardware in finding reliable solutions and the extent to which quantum annealers represent (or have the potential to represent) a realistic alternative to simulated annealing and deterministic approaches in the context of this type of combinatorial problems. We also emphasize that investigating the scalability of quantum applications using moderate-to-small scale test problems is particularly important in the context of quantum-classical hybrid approaches based on the quantum annealing hardware. Particularly, the quantum component can be used as a specialized sub-processor for improving sub-problems within larger classical algorithms \cite{hybrid2, optimizingtheoptimizer, hybrid}.

\section{Related work}
A variety of work has been done in the area of drone scheduling, as well as in related optimization problems with similar formulations. Examples, concerned with the fields of transportation, logistics, and supply chain management, include applications of the renowned Travelling Salesperson Problem (TSP), for which decomposition techniques are introduced in e.g.\ \cite{veh2, veh4, veh}. Different versions have been analyzed, both for minimizing emissions \cite{veh}, and addressing the dynamic version with real-time customer requests \cite{veh2}.

Another widely studied problem is the Bin Packing Problem (BPP), where items of different sizes must be allocated into a finite number of bins with a fixed given capacity, so to minimize bin usage. For this scope, a hybrid classical-quantum approach is presented in \cite{bpp}. The parallelism with the DDPP is clear: in the DDPP case, the challenge is efficiently utilizing drones to deliver a set of parcels while considering battery capacities; similarly, in the classic BPP, the aim is to efficiently pack a set of items into bins, while considering limited physical capacities. A special application \cite{bppfuel} is the optimization of spent nuclear fuel (SNF) filling in canisters, so that the maximum heat output does not exceed a limiting value. Another scheduling-related problem is the optimal flight gate assignment for airport management, aiming to the minimization of the total transit time of the passengers in an airport, with time constraints \cite{traffic, scheduling, gate}.

The starting points for our work are the results of Jana and Mandal in \cite{bppjanamandal}, where they prove the NP-hardness of the DDPP, and give an ILP representation for which they propose two approximation algorithms. Specifically, they give a greedy approximation algorithm, and a colouring-based algorithm, on a graph representation of the DDPP. Indeed, they construct, for the given set of delivery time intervals, an interval graph, where vertices represent intervals and are adjacent if the corresponding intervals conflict.

In our work, we use their ILP formulation as a starting point to develop two QUBO formulations of the DDPP, amenable to global optimization via sampling for both classical and quantum annealing. 

The following section is dedicated to mathematically formalizing the DDPP, as a general constrained optimization problem first, then as an integer linear programming problem, and finally as a QUBO problem. The section presents two QUBO formulations: one obtained following a relatively standard variable augmentation approach \cite{veh2, qubo2} and another one obtained via a tailored relaxation. An analysis is provided showing that the latter formulation improves the former in terms of the number of variables and qubits required, positively impacting the time-to-solution of annealing optimization samplers. 

\section{QUBO formulation of the drone delivery packing problem}\label{qubomodels}
In this section, we introduce the drone delivery packing problem as a combinatorial optimization problem and we then derive two alternative QUBO formulations.

The DDPP seeks to give an optimal delivery assignment to a set of identical available drones, with battery budget. Their task is to bring to completion a certain set of deliveries requested by customers, with constraints regarding cost and time consistency.

More specifically, for the given set of deliveries, their delivery time intervals, cost in terms of energy for each delivery, and
(equal) battery budget of the drones, the goal is to schedule an optimal set of drones, that minimizes the used drones.

The assignment must guarantee time and energy consistency, other than the completion of the deliveries. Let $\mathcal{M} = \{1,...,m\}$ be the collection of identical drones, available for the company in a given depot. Let $\mathcal{N} = \{1, ..., N\}$ be the set of deliveries assigned to the company to be completed. The binary variables used within the model are: 
\begin{align*}
    &x_{ij} = \begin{cases}
                1, \ \text{drone $i \in \mathcal{M}$ delivers to $j \in \mathcal{N}$}  \vspace{0.2mm} \\
                0, \ \text{otherwise}
        \end{cases} 
\end{align*}
Other variables involved in the constraints formalization are: 
    \begin{enumerate}
        \item The battery budget $B > 0$ for each drone, namely the service duration time of the battery;
        \item The delivery time intervals $I_j$, representing the time window in which delivery $j \in \mathcal{N}$ is done; 
        \item The cost $c_j$ in terms of battery for completing delivery $j \in \mathcal{N}$.
    \end{enumerate}
Then, the DDPP can be formulated as follows \cite{bppjanamandal}: 
\begin{align}\label{eq:ILP1}  
&\min_{x_{ij} \in \{0,1\}} \sum_{i \in \mathcal{M}} \max_{j\in \mathcal N} x_{ij}\, , \quad \text{s.t.} \\
&\begin{cases}
    \sum_{j \in \mathcal{N}} c_j x_{ij} \le B & \forall \ i \in \mathcal{M}\\
    \sum_{i \in \mathcal{M}} x_{ij} = 1 & \forall \ j \in \mathcal{N}  \\
    x_{ij} + x_{ik} \le 1 & \forall \ i \in \mathcal{M}, \\ 
    & j, k \in \mathcal{N}, I_j \cap I_k \neq \emptyset 
\end{cases}\notag
\end{align}
Notice that the binary quantity $\max_ {j\in \mathcal N} x_{ij}$ indicates whether drone $i$ is used. In fact,
\[
\max_{j\in \mathcal N} x_{ij}= \begin{cases}
    1 & \text{if drone $i$ has been used} \\
    0, & \text{otherwise}
\end{cases}
\]
Therefore, $\sum_{i \in \mathcal{M}} \max_{j\in \mathcal N} x_{ij}$ counts the number of used drones.

\subsection{Standard ILP form of the maximum function and first QUBO formulation}

Starting from the combinatorial optimization formulation \eqref{eq:ILP1}, we are interested in obtaining an equivalent QUBO formulation representing the Hamiltonian to be given to the annealer. This can be done with two relatively standard steps: (see e.g.\ \cite{veh2, veh, qubo2, qubo1}):
\begin{enumerate}
\item Transforming all the inequality constraints into equality constraints by adding a number of slack variables, thus reformulating \eqref{eq:ILP1} as an integer (binary, to be precise) linear programming (ILP).
\item From the binary ILP obtain the QUBO by transferring all the linear constraints into penalty terms to be added to the ILP objective function (which is quadratic and binary).
\end{enumerate}
Note that, for each inequality constraint, the number of slack variables comes from the quantity $s = b - \sum_{i=1}^N A_ix_i \le b - \sum_{i : A_i \le 0} A_i$. Therefore, at most $\lceil \log_2 \big(b - \sum_{i : A_i \le 0} A_i \big) \rceil$ coefficients are needed to represent $s$, and each inequality constraint requires this number of slack variables for the transition to equality.

In order to apply steps 1 and 2 above to the DDPP model \eqref{eq:ILP1}, we first notice that, while the constraints in \Cref{eq:ILP1} are linear inequalities, the objective function is written in terms of $\max_j \sum_i x_{ij}$ and thus it is not quadratic. 
In order to formulate the problem using a quadratic binary loss, we introduce $m$ new binary variables $y_i$ and we notice that the condition $y_i := \max_j \sum_i x_{ij}$ (for binary $x_{ij}$ and $y_i$) corresponds to the pair of conditions $x_{ij}\leq y_j \qquad \text{and}\qquad y_i \leq \sum_{j\in \mathcal N} x_{ij}$ for all $i\in \mathcal M$ and all $j\in\mathcal N$. By inserting these constraints into \eqref{eq:ILP1} we obtain a new constrained quadratic binary problem
\begin{align}
&\min_{y_i \in \{0,1\}} \sum_{i \in \mathcal{M}} y_i\, , \quad \text{s.t.}\label{eq:ILP2} \\
&\begin{cases}
\sum_{j \in \mathcal{N}} c_j x_{ij} \le B & \forall \ i \in \mathcal{M} \\
\sum_{i \in \mathcal{M}} x_{ij} = 1 & \forall \ j \in \mathcal{N} \\
x_{ij} + x_{ik} \le 1 & \forall \ i \in \mathcal{M}, \\ 
    & j, k \in \mathcal{N}, I_j \cap I_k \neq \emptyset  \\
x_{ij} \le y_i & \forall \ i \in \mathcal{M}, \forall \ j \in \mathcal{N} \\
y_i \le \sum_{j\in \mathcal N} x_{ij} & \forall \ i \in \mathcal{M}
\end{cases}\notag
\end{align}

Now, we transform the inequality constraints in \eqref{eq:ILP2} into quadratic equality constraints. This is done by means of a number of additional slack variables, obtaining the following quadratic penalty terms \cite{penalties}:
\begin{equation*}
    \begin{aligned}
        &H_{C_1} =  \sum_{i \in \mathcal{M}} \bigg( \sum_{j \in \mathcal{N}} c_j x_{ij} + \sum_{l=0}^{\lceil \log_2 B \rceil - 1} 2^l s_{il} - B \bigg)^2 \\
        &H_{C_2} = \sum_{j \in \mathcal{N}} \bigg(  \sum_{i \in \mathcal{M}} x_{ij} - 1 \bigg)^2 \\
        &H_{C_3} =  \sum_{i \in \mathcal{M}} \sum_{j,k \in \mathcal{N} : I_j \cap I_k \neq \emptyset} \bigg( x_{ij} + x_{ik} + t_{i} - 1 \bigg)^2 \\
        &H_{C_4} =   \sum_{i \in \mathcal{M}} \sum_{j \in \mathcal{N}} \bigg( x_{ij} - y_i + r_{ij} \bigg)^2 \\
        &H_{C_5} =  \sum_{i \in \mathcal{M}} \bigg( - \sum_{j \in \mathcal{N}} x_{ij} + y_i + \sum_{l=0}^{\lceil \log_2 ( |\mathcal{N}|) \rceil -1} 2^l p_{il}\bigg)^2 \\
    \end{aligned}
\end{equation*}
Finally, by moving the binary constraints $H_{C_i}$ into the objective function as penalty terms we obtain the  Hamiltonian
\begin{equation}    \label{eq:QUBO1}
       H = \sum_{i \in \mathcal{M}} y_i  + \sum_{i=1}^5 \alpha_i H_{C_i} 
\end{equation}
where $\alpha_i > 0$, for  $i=1,2,3,4,5$. This Hamiltonian, to be given to the annealer, is made of six terms, the objective function and the five quadratic Hamiltonian terms $H_{C_i}$, and is a function of the original variables $x_{ij}$ and several additional slack variables $y_i$, $c_j$, $t_i$, $r_{ij}$, and $p_{ij}$. 
With this construction, the overall number of variables needed to formulate DDPP as a QUBO problem is 
\begin{equation}
    \label{eq:num_var}
    \begin{aligned}
&|\mathcal{M}| + |\mathcal{M}||\mathcal{N}| + |\mathcal{M}| \lceil \log_2 B \rceil  + |\mathcal{M}|  \kappa + \\ &|\mathcal{M}||\mathcal{N}| + |\mathcal M| \lceil \log_2 (|\mathcal{N}|) \rceil 
    \end{aligned}
\end{equation}
where $\kappa$ denotes the number of time-conflicting deliveries: $ \kappa = |\bigcup_{j,k \in \mathcal{N}} \{j, k : i < j, I_i \cap I_j \neq 0\}|$.

The first two terms come from the original sets of variables $\{y_i\}_{i \in \mathcal{M}}$ and $\{x_{ij}\}_{i \in \mathcal{M}, j \in \mathcal{N}}$, respectively, while the subsequent terms come from $H_{C_1}$, $H_{C_3}$, $H_{C_4}$ and $H_{C_5}$. While the obtained QUBO formulation results from the most standard method, we notice that it leads to a number of overall variables that quickly grow very large.

In the next subsection, we present an alternative QUBO model, that better leverages the structure of the problem, while significantly reducing the number of variables.  

\subsection{Quadratic proxy of the maximum function and second QUBO formulation}\label{QUBOformulation}

\begin{figure}[t]
    \begin{tabular}{@{}lcrc@{}}
    {\!\!\!\!\LARGE $[x_{ij}]$} & {\LARGE $\gets$} &
        $\begin{bmatrix}
            0 & 1 & 0 & \cdots \\
            1 & 0 & 0 & \cdots \\
            0 & 0 & 1 & \cdots \\
            \vdots & \vdots & \vdots & \ddots 
        \end{bmatrix}$ & $\begin{bmatrix}
            1 & 1 & 1 & \cdots \\
            0 & 0 & 0 & \cdots \\
            0 & 0 & 0 & \cdots \\
            \vdots & \vdots & \vdots & \ddots 
        \end{bmatrix} $
    \end{tabular}
    \caption{Example of worst (left) and best (right) configurations of $x_{ij}$ for \eqref{eq:ILP1} if only the $\sum_i x_{ij}=1$ constraint were to be considered. }
    \label{fig:x_ij}
\end{figure}

In this subsection, we propose a more efficient QUBO model for the DDPP that differs from the one introduced in the previous section. The key idea is to consider a cheaper quadratic objective function that acts as a proxy for the original $\max_j x_{ij}$ without requiring additional slack terms, thus directly expressing the optimization problem in terms of the initial $x_{ij}$ variables. 

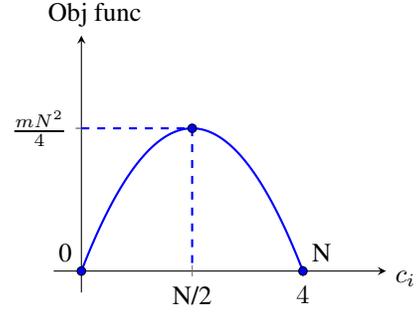
\begin{figure}[t]
\centering
\begin{tikzpicture}
\begin{axis}[
    height=5cm,
    width=6cm,
    axis lines=middle,
    xlabel={$c_i$},
    ylabel={Obj func},
    ylabel near ticks,
    xlabel near ticks,
    xtick={0, 2, 4},
    ytick={0, 4},
    yticklabels={0, $\frac{mN^2}{4}$},
    xmin=0, xmax=5,
    ymin=0, ymax=6,
    enlargelimits=true,
    domain=0:4,
    samples=100,
    ylabel style={at={(ticklabel* cs:1)}, anchor=south, rotate=-90, xshift=5pt},
    xlabel style={at={(ticklabel* cs:1)}, anchor=west, yshift=-3pt},
    extra x ticks={2},
    extra x tick style={grid=none, tick label style={below, fill=white}},
    extra x tick labels={N/2},
    extra y ticks={4},
    extra y tick style={grid=none, tick label style={left, fill=white}},
    extra y tick labels={$\frac{mN^2}{4}$},
] 
    \addplot[
        only marks,
        mark=*,
        mark options={scale=0.8, fill=blue}
        ]
        coordinates {
        (0,0)
        (2,4)
        (4,0)
        };
    \addplot[blue, thick, dashed] coordinates {(2,4) (2,0)};
    \addplot[blue, thick, dashed] coordinates {(0,4) (2,4)};
    \addplot[blue, thick] {(4-x)*x};
    \node[anchor=south east] at (axis cs:0,0) {0};
    \node[anchor=south west] at (axis cs:4,0) {N};
\end{axis}
\end{tikzpicture}
\caption{Continuous approximation in dimension 2 of the parabola $\sum_{i \in \mathcal{M}} \Big(N- \sum_{j\in \mathcal M}x_{ij}\Big)\sum_{j\in \mathcal M}x_{ij}$ as a function $\mathbb{N}^m \to \mathbb{N}$ of the $m$ integer values $c_i = \sum_{j\in \mathcal M}x_{ij}$ describing for $i \in \{1,...,m\}$ the number of deliveries completed by drone $i$.}
\label{fig:parabola}
\end{figure}

To this end, we notice that in terms of the assignment matrix $[x_{ij}]$, seeking to minimize $\sum_i \max_j x_{ij}$ as in the original minimization problem \eqref{eq:ILP1} is equivalent to favour configurations of $[x_{ij}]$ with as many zero rows as possible. In other words, taking into account the ``all deliveries once'' condition $\sum_i x_{ij}=1$ and ignoring the other constraints, the worst configuration for $[x_{ij}]$ is one with a permutation-like pattern, while the best configuration is one where all rows in $[x_{ij}]$ are zero except for one, i.e.\ all the deliveries are done by only one drone, as illustrated in \Cref{fig:x_ij}. If $\sigma_i = \sum_j x_{ij}$ denotes the row-sum of $[x_{ij}]$, we then notice that the quadratic function $\sum_i (N-\sigma_i)\sigma_i$ combined with the all-delivery-once constraints favours the latter type of configuration and penalizes the former, as $(N-\sigma_i)\sigma_i$ is minimal when either $\sigma_i=0$ or $\sigma_i=N$. See \Cref{fig:parabola} for an illustration.

With these considerations in mind, we propose here the following modified QUBO formulation for the DDPP
\begin{align}   
&\min_{x_{ij} \in \{0,1\}} \sum_{i \in \mathcal{M}} \Big(N- \sum_{j\in \mathcal M}x_{ij}\Big)\sum_{j\in \mathcal M}x_{ij}\, , \quad \text{s.t.}\label{eq:new-ILP}\\
&\begin{cases}
    \sum_{j \in \mathcal{N}} c_j x_{ij} \le B & \forall i \in \mathcal{M}\\
    \sum_{i \in \mathcal{M}} x_{ij} = 1 & \forall j \in \mathcal{N}  \\
    x_{ij} + x_{ik} \le 1 & \forall i \in \mathcal{M}, \ j, k \in \mathcal{N}, I_j \cap I_k \neq \emptyset 
\end{cases}\notag
\end{align}
which can be expressed in QUBO format by means of the quadratic unconstrained Hamiltonian
\begin{equation}\label{eq:QUBO2}
    H = H_0+\alpha_1 H_{C_1}+\alpha_2H_{C_2}+\alpha_3 H_{C_3},
\end{equation}
with $H_0 = \sum_{i} (N- \sum_{j}x_{ij})\sum_{j}x_{ij}$.

Note that this new model does not need to consider the variables $y_i$, nor the consistency constraints. Thus, it requires only
\begin{equation}
    \label{eq:num_var_improved}
    |\mathcal{M}||\mathcal{N}| + |\mathcal{M}| \lceil \log_2 B \rceil  + |\mathcal{M}|  \kappa
\end{equation}
variables, as opposed to \eqref{eq:num_var}. 

\begin{figure}[t]
\begin{tikzpicture}
\begin{axis}[
    xlabel={$N$ (number of deliveries)},
    ylabel={Num variables},
    ylabel near ticks,
    xlabel near ticks,
    xmin=4, xmax=8,
    ymin=50, ymax=160,
    xtick={4,5,6,7,8},
    ytick distance=30,
    grid=both,
    axis lines=middle,
    enlargelimits=true,
    title={},
    legend style={at={(0.5,1.05)},anchor=south},
    ylabel style={xshift=-20pt, yshift=15pt},
    xlabel style={yshift=-32pt, xshift=-35pt} 
    ]

\addplot[
    only marks,
    mark=*,
    mark options={scale=0.8, fill=pink}
] 
coordinates {
    (4,55)
    (5,60)
    (6,65)
    (7,70)
    (8,75)
};

\addplot[
    only marks,
    mark=*,
    mark options={scale=0.8, fill=violet}
] 
coordinates {
    (4,85)
    (5,100)
    (6,110)
    (7,120)
    (8,130)
};

\addlegendentry{QUBO \eqref{eq:QUBO2}}
\addlegendentry{QUBO \eqref{eq:QUBO1}}

\addplot[
    pink,
    thick,
    mark=none
] 
coordinates {
    (4,55)
    (5,60)
    (6,65)
    (7,70)
    (8,75)
};

\addplot[
    violet,
    thick,
    mark=none
] 
coordinates {
    (4,85)
    (5,100)
    (6,110)
    (7,120)
    (8,130)
};

\end{axis}
\end{tikzpicture}
\caption{Comparison between the number of logical variables necessary to encode the model \eqref{eq:QUBO1} and \eqref{eq:QUBO2}. 
We let $N$ vary in $\{4,\dots,8\}$, we chose a fixed $m=5$ and a fixed battery budget $B = 5$. Costs $c_j$ for each delivery $j \in N$ are randomly sampled from a Gaussian distribution such that  $0\leq c_j \le B \ \forall \ j \in N$. Time intervals $I_j$ are sampled from a uniform distribution from $8$ a.m. to $8$ p.m. and are at least one hour long and at most two.}
\label{fig:qubits}
\end{figure}

The advantage of the modified QUBO formulation is further highlighted in \Cref{fig:qubits} in terms of logical variables involved within the formulation. Two curves are shown, quantifying the improvement in the number of logical variables when using \Cref{eq:QUBO2} with respect to \Cref{eq:QUBO1}, on example problems with $N\in \{4,5,6,7,8\}$. This also shows how the benefit in the variables usage scales with increasing deliveries. In particular, we notice that the variable allocation decreases significantly and that the discrepancy between the two models grows as $N$ increases.

The following section moves on to Simulated and Quantum Annealing techniques for tackling the DDPP formulated as \Cref{eq:QUBO2}, while addressing the consequent state-of-the-art quantum technology limitations.

\section{Simulated and quantum annealing}\label{qaperf}

In this section, we start by introducing the main ideas at the basis of the annealing optimization algorithm, highlighting the differences between the classical and the quantum versions and, in particular, the potential advantages of the latter over the former. 

Simulated and Quantum Annealing are probabilistic algorithms for approximating minimal solutions by exploring the loss landscape and the solution space iteratively, through certain acceptance probabilities. Inspired by the physical phenomenon of metallurgical annealing, the process leverages the adiabatic theorem (\cite{adiab}). 

Therefore, the system is initialized to an easily prepared state and slowly evolves, through a time-dependent perturbation, towards a more complex Hamiltonian, whose ground state encodes the solution to the optimization problem. More specifically, the algorithm aims at finding the global minimum of an objective function. It starts with a guess and iteratively moves to neighboring states. The switch is based on a Boltzmann acceptance probability that exponentially depends on the height of the energy hill representing the local minimum to escape.

In the Quantum version, the perturbation is a time-dependent magnetic field, allowing quantum phenomena like the tunneling effect (see \cite{physics}). This modifies the acceptance probability that now depends also on the length of the potential barrier.

The extra dependence of QA's acceptance probability on $L$ allows the system to compensate for the difficulty of escaping from very high barriers, e.g. very isolated local minima. This shows that QA can potentially outperform SA, especially when the loss landscape consists of very high but thin barriers surrounding shallow local minima. Nevertheless, this is problem-specific, as the acceptance probability is strongly related to the shape of the objective function.

Rigorous proofs of the advantage of QA over SA in terms of the ability to escape local optima exist only for specific examples \cite{qaoutperformance2, qaoutperformance4, qaoutperformance}. However, there is a wealth of empirical indications that quantum annealing can make highly complex combinatorial optimization problems computationally feasible, see e.g. \cite{Boixo, qaenhancement, qaenhancement2, Kadowaki_1998, higham2022testing,bergermann2024nonlinear}. Also, transferring this theoretical advantage of QA over SA in terms of time-to-solution is a challenging task. In fact, the relationship between acceptance probability and annealing time is still an open question \cite{nature_nonad, vanilla}. 

However, we will show in the following subsection the time-to-solutions found for our example problems both for SA and QA, displaying a significant gain for the considered problem sizes. Also, the next subsection starts by addressing the issue of embedding a QUBO in the chipset of a Quantum Annealer both in general and for the QUBO formulation \eqref{eq:QUBO2}.

\subsection{Chipset embedding}\label{emb}

In this subsection, we provide a brief overview of the quantum annealer hardware, and its connection with the QUBO format.
The QUBO formulation is strongly related to the quantum annealing computation. In fact, the processor of a quantum annealer can be seen as a graph whose nodes are the qubits and whose edges represent the connections between the qubits. Similarly, QUBO formulations, being quadratic forms, can be represented as a graph where linear and quadratic terms are associated with the nodes and the edge weights of the graph, respectively. 

Thus, QUBO problems are particularly suitable to be implemented and solved on quantum annealers, as long as an embedding of the QUBO graph into a subgraph of the fixed quantum annealer's chipset graph exists. However, QPU chipsets have specific connectivity limitations due to the physical layout of their qubits and the coupling strengths between them. Thus, within the embedding process, it is often the case that multiple physical qubits are involved to represent single model variables, forming so-called qubit chains (\cite{chainstrength}). This leads to the concrete involvement of a greater number of physical qubits, even though, at the model level, fewer variables were used. This gap between model variables and the physical qubits involved can be significant. Nevertheless, it is intended to decrease with the technological advancement, and the construction of bigger and more connected quantum computers, allowing more direct embeddings. We remark that the D-Wave Leap's Ocean SDK contains several specifically designed utilities for this task. See for example \cite{minors} or D-Wave's documentation  \cite{documentation}. 

\begin{figure}[t]
\begin{tikzpicture}
\begin{axis}[
    xlabel={$N$ (number of deliveries)},
    ylabel={Num qubits},
    ylabel near ticks,
    xlabel near ticks,
    xmin=4, xmax=8,
    ymin=100, ymax=440,
    xtick={4,5,6,7,8},
    ytick distance=30,
    grid=both,
    axis lines=middle,
    enlargelimits=true,
    title={},
    legend style={at={(0.5,1.05)},anchor=south},
    ylabel style={xshift=-20pt, yshift=15pt},
    xlabel style={yshift=-32pt, xshift=-35pt} 
    ]

\addplot[
    only marks,
    mark=square*,
    mark options={scale=0.8, fill=pink}
] 
coordinates {
    (4,129)
    (5,188)
    (6,228)
    (7,253)
    (8,316)
};

\addplot[
    only marks,
    mark=square*,
    mark options={scale=0.8, fill=violet}
] 
coordinates {
    (4,184)
    (5,275)
    (6,309)
    (7,358)
    (8,434)
};

\addlegendentry{QUBO \eqref{eq:QUBO2}}
\addlegendentry{QUBO \eqref{eq:QUBO1}}

\addplot[
    pink,
    thick,
    mark=none
] 
coordinates {
    (4,129)
    (5,188)
    (6,228)
    (7,253)
    (8,316)
};

\addplot[
    violet,
    thick,
    mark=none
] 
coordinates {
    (4,184)
    (5,275)
    (6,309)
    (7,358)
    (8,434)
};

\end{axis}
\end{tikzpicture}
\caption{Comparison between average number of physical qubits involved when using models \eqref{eq:QUBO1} and \eqref{eq:QUBO2}. We let $N$ vary in $\{4,\dots,8\}$, we chose a fixed $m=5$ and a fixed battery budget $B = 5$. Costs $c_j$ for each delivery $j \in N$ are randomly sampled from a Gaussian distribution such that $0 \leq c_j \leq B \ \forall \ j \in N$. Time intervals $I_j$ are sampled from a uniform distribution from 8 a.m. to 8 p.m. and are at least one hour long and at most two.}
\label{fig:physical_qubits}
\end{figure}
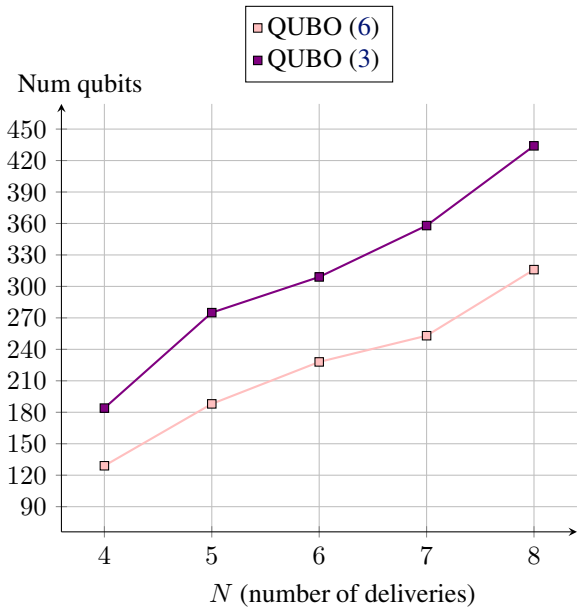

Also, we point out that the advantage obtained when using \Cref{eq:QUBO2} over \Cref{eq:QUBO1} discussed in \Cref{QUBOformulation}, results in a reduction in the number of variables to be embedded in the chosen quantum hardware. This highlights how choosing a more convenient model formulation extends to a meaningful decrease in the number of variables involved, which is even more evident at a hardware level. This is shown in \Cref{fig:physical_qubits}, where a "hardware counterpart" of \Cref{fig:qubits} is depicted. In fact, while the former shows an improvement in terms of logical variables, the latter shows similar plots of variables, but referring to physical ones. Also, this reduction is observed to scale with the number $N$ of deliveries, as in \Cref{fig:qubits}. Note that the qubits usage depends on the chosen embedding. We used the EmbeddingComposite class provided by D-Wave's Ocean SDK for this process \cite{documentation, emb}. This tool automatically minor-embeds the problem into the device, calculating a new minor-embedding each time one of its sampling methods is called \cite{minors}. Therefore, the values shown in \Cref{fig:physical_qubits} represent average numbers of physical qubits, computed across 10 different embedding computations. In this sense, such values measure how many physical variables are needed, on average, to encode and embed each problem instance for increasing $N \in \{4,...,8\}$.

The following section is dedicated to introducing the experimental configuration used for conducting our tests. Particularly, we describe the parameter setting, the hardware specifications both of the local processor used for running simulated annealing and the quantum processor, and the synthetic problem instances generated for the testing.

\section{Numerical setup and parameter tuning}\label{sec:experiments}

\subsection{Choice of the penalty coefficients}\label{penalties}

\begin{table*}[t]
\centering
\tiny
\setlength\tabcolsep{2pt}
\begin{tabularx}{\textwidth}{|l|l|l|l|X|X|}
\hline
\multicolumn{1}{|c|}{\textbf{Instance}} & \multicolumn{1}{c|}{\textbf{\# deliveries}} & \multicolumn{1}{c|}{\textbf{battery budget}} & \multicolumn{1}{c|}{\textbf{distribution}} & \multicolumn{1}{c|}{\textbf{costs}}                                          & \multicolumn{1}{c|}{\textbf{time intervals}}                                                                                                                                 \\ \hline
1                                       & 10                                          & 70                                           & gaussian                                   & {[}59.8, 42.2, 27.4, 40.3, 43.6, 30.0, 31.7, 41.7, 9.6, 34.1{]}              & {[}{[}14, 16{]}, {[}9, 10{]}, {[}14, 17{]}, {[}16, 18{]}, {[}13, 14{]}, {[}10, 11{]}, {[}13, 16{]}, {[}13, 16{]}, {[}11, 12{]}, {[}8, 10{]}{]}                               \\ \hline
2                                       & 10                                          & 100                                          & uniform                                    & {[}7.0, 31.5, 44.8, 90.5, 9.3, 14.2, 79.0, 2.1, 91.1, 57.3{]}                & {[}{[}12, 14{]}, {[}8, 10{]}, {[}13, 16{]}, {[}15, 16{]}, {[}17, 20{]}, {[}16, 19{]}, {[}13, 14{]}, {[}14, 17{]}, {[}16, 18{]}, {[}18, 20{]}{]}                              \\ \hline
3                                       & 12                                          & 70                                           & gaussian                                   & {[}29.8, 24.9, 45.5, 32.3, 23.5, 22.2, 42.1, 35.5, 48.6, 54.2, 39.1, 33.1{]} & {[}{[}18, 19{]}, {[}19, 20{]}, {[}15, 18{]}, {[}12, 13{]}, {[}13, 15{]}, {[}16, 18{]}, {[}9, 12{]}, {[}15, 17{]}, {[}17, 19{]}, {[}8, 9{]}, {[}17, 18{]}, {[}18, 20{]}{]}    \\ \hline
4                                       & 10                                          & 70                                           & gaussian                                   & {[}42.0, 54.9, 32.3, 39.7, 13.9, 9.0, 21.3, 35.5, 48.9, 33.4{]}              & {[}{[}14, 17{]}, {[}16, 19{]}, {[}8, 11{]}, {[}12, 14{]}, {[}18, 19{]}, {[}16, 18{]}, {[}10, 11{]}, {[}16, 18{]}, {[}14, 16{]}, {[}8, 9{]}{]}                                \\ \hline
5                                       & 10                                          & 70                                           & gaussian                                   & {[}34.0, 30.9, 31.1, 38.7, 27.1, 20.6, 29.5, 44.2, 38.3, 26.4{]}             & {[}{[}8, 10{]}, {[}13, 16{]}, {[}12, 14{]}, {[}18, 20{]}, {[}8, 9{]}, {[}17, 18{]}, {[}16, 18{]}, {[}19, 20{]}, {[}9, 11{]}, {[}8, 11{]}{]}                                  \\ \hline
6                                       & 12                                          & 50                                           & uniform                                    & {[}1.5, 16.1, 47.8, 1.7, 28.1, 20.7, 1.2, 31.8, 8.2, 48.2, 3.2, 43.5{]}      & {[}{[}10, 12{]}, {[}18, 19{]}, {[}19, 20{]}, {[}18, 19{]}, {[}16, 18{]}, {[}16, 18{]}, {[}10, 13{]}, {[}12, 14{]}, {[}13, 14{]}, {[}16, 18{]}, {[}16, 17{]}, {[}18, 19{]}{]} \\ \hline
7                                       & 12                                          & 70                                           & gaussian                                   & {[}28.2, 31.6, 44.5, 35.5, 41.3, 37.0, 20.7, 47.1, 33.3, 29.7, 40.8, 25.3{]} & {[}{[}19, 20{]}, {[}15, 18{]}, {[}8, 10{]}, {[}15, 18{]}, {[}19, 20{]}, {[}12, 14{]}, {[}8, 9{]}, {[}16, 18{]}, {[}16, 19{]}, {[}16, 17{]}, {[}13, 16{]}, {[}9, 12{]}{]}     \\ \hline
8                                       & 12                                          & 50                                           & uniform                                    & {[}4.6, 46.5, 12.2, 44.3, 13.2, 2.2, 25.0, 47.8, 25.3, 37.9, 44.6, 0.9{]}    & {[}{[}12, 15{]}, {[}11, 12{]}, {[}17, 19{]}, {[}9, 11{]}, {[}16, 19{]}, {[}10, 13{]}, {[}10, 12{]}, {[}15, 18{]}, {[}14, 15{]}, {[}12, 15{]}, {[}17, 19{]}, {[}13, 15{]}{]}  \\ \hline
9                                       & 10                                          & 100                                          & gaussian                                   & {[}43.9, 47.9, 55.2, 47.1, 52.6, 53.7, 54.0, 48.3, 41.1, 46.7{]}             & {[}{[}17, 20{]}, {[}8, 11{]}, {[}17, 19{]}, {[}13, 16{]}, {[}17, 18{]}, {[}19, 20{]}, {[}13, 15{]}, {[}13, 16{]}, {[}14, 16{]}, {[}13, 14{]}{]}                              \\ \hline
10                                      & 12                                          & 70                                           & gaussian                                   & {[}30.6, 18.1, 27.8, 18.7, 40.2, 34.3, 39.6, 31.1, 41.5, 40.9, 41.6, 27.5{]} & {[}{[}10, 12{]}, {[}8, 10{]}, {[}14, 17{]}, {[}12, 13{]}, {[}16, 18{]}, {[}13, 14{]}, {[}12, 13{]}, {[}12, 15{]}, {[}17, 19{]}, {[}8, 10{]}, {[}8, 11{]}, {[}14, 16{]}{]}    \\ \hline
11                                      & 10                                          & 70                                           & uniform                                    & {[}17.0, 31.2, 58.9, 39.4, 24.1, 66.4, 31.9, 40.2, 64.6, 23.6{]}             & {[}{[}8, 9{]}, {[}9, 10{]}, {[}14, 17{]}, {[}15, 18{]}, {[}10, 12{]}, {[}16, 17{]}, {[}15, 18{]}, {[}10, 12{]}, {[}13, 14{]}, {[}17, 18{]}{]}                                \\ \hline
12                                      & 12                                          & 100                                          & uniform                                    & {[}70.6, 90.1, 1.4, 73.0, 22.8, 88.8, 87.2, 28.5, 99.5, 35.8, 4.3, 1.3{]}    & {[}{[}18, 19{]}, {[}11, 12{]}, {[}12, 15{]}, {[}10, 13{]}, {[}16, 17{]}, {[}12, 14{]}, {[}15, 16{]}, {[}12, 15{]}, {[}9, 10{]}, {[}15, 18{]}, {[}16, 19{]}, {[}12, 13{]}{]}  \\ \hline
\end{tabularx}
\caption{Sampled problem sets. Drones are fixed to be $10$, deliveries $N$ are randomly picked in $\{10,12\}$ and battery budgets in $\{50, 70, 100\}$. Also, a distribution is randomly chosen between uniform and Gaussian for each instance, and according to this, we generate $N$-dimensional lists of delivery energy costs and delivery time intervals, respectively. Time windows are generated from 8 a.m. to 8 p.m. and are at least one and at most two hours long, while costs $c_j$ are required to be smaller than the budget.}
\label{tab:k90_table}
\end{table*}

In this subsection, we focus on the choice of the penalty coefficients $\alpha_i$ in \eqref{eq:QUBO2}. We recall that $\alpha_1$ is responsible for the energy budget constraint $H_{C_1}$; $\alpha_2$ is related to the requirement that all deliveries should be done once, as modeled by the constraint term $H_{C_2}$; $\alpha_3$ is responsible for the time interval constraints as in $H_{C_3}$.
In general, higher penalty coefficients emphasize constraint satisfaction, at the expense of sacrificing optimality, while lower coefficients may prioritize the objective function over feasibility. For this reason, the process of setting penalty weights is not trivial, and it becomes more challenging when dealing with multiple constraints. In fact, in that case, one has to choose the coefficients to balance the importance of individual terms, and that is typically problem-specific. 

A common approach to address the interplay between single penalty terms (see e.g.\ \cite{penalties}) is to set all the weights $\alpha_i$ to a value greater than the largest possible absolute value of the original objective function $H_{0}$. In this way, violating constraints is considered more significant than minimizing the objective function. In the case of \Cref{eq:QUBO2}, the objective function $H_0$ represents a discrete and positive $m$-dimensional downward-facing parabola touching the $x$-axis in $0$ and $N$, as depicted in \Cref{fig:parabola}. Therefore, its maximum value is upper bounded as $H_0 \leq m \frac{N^2}{4}$.

In our case, however, setting the penalty coefficients uniformly to the same value did not work well as the resulting solutions consistently failed to satisfy the constraint imposing that all the deliveries should be performed by exactly one drone. For this reason, we performed a parameter tuning analysis for $\alpha_2$ on a range of randomly generated benchmark test problems. We fixed the number of drones to $m=10$ and sampled problem instances with either $10$ or $20$ deliveries, $N\in \{10,12\}$, by randomly generating costs and time intervals. \Cref{tab:k90_table} shows the details of the test problems we considered, and more details on the generation process can be found in the next  \Cref{sec:optimization-performance}. In order to tune the parameter $\alpha_2$ we set
\begin{equation}\label{eq:constraint_all-at-once}
    \alpha_2 = k \cdot (mN^2)/4
\end{equation}
and let $k$ vary within a discrete grid from $5$ to $150$. We tested the performance of SA on each of these parameter values for the dataset from \Cref{tab:k90_table} and checked the quality of the solution in terms of objective function and constraint satisfaction. As a result, we observed overall best performance within the considered problem instances dimension range for $k=120$. 

The next subsection provides additional details on the set of test problems as well as an in-depth analysis of the performance of both quantum and simulated annealing samplers, as compared to the exact solution computed with a commercial deterministic solver (Gurobi \cite{gurobi}, used with an academic license) applied to the original constrained ILP \eqref{eq:ILP1}. In all our tests, we fixed the value of $\alpha_2$ using $k=120$ and set $\alpha_i=mN^2/4$ for the other coefficients, $i \neq 2$. 

For what concerns the Quantum Annealer, we conducted an independent tuning of penalty weights, to empirically find their optimal combination, also for this different search method. In fact, the coefficients setting requires a different selection, due to the different exploration of the solution space and the different strategies used to escape local minima. We discovered that also QA encounters more difficulties in fulfilling the all-deliveries-once constraint and that the same strengthening coefficient \eqref{eq:constraint_all-at-once} with $k=120$ performs well in this case.

\subsection{Test problems setup}\label{sec:problems-setup}

In order to generate an unbiased test problem benchmark, we use a random instance generator to sample 24 synthetically generated problem settings encompassing a range of possible configurations. This is a standard benchmarking approach, see e.g.\ \cite{bpp}. In particular, we produce two datasets, with the aim of analyzing the performance sensitivity on the problem dimension, especially with respect to the deliveries set size. The first dataset comprises 12 instances with either 10 or 12 deliveries, resulting in a larger number of variables, while the second group of 12 instances uses either 7 or 8 deliveries, leading to fewer variables. We report in \Cref{tab:k90_table} the sampled sets for the cases of 10 and 12 deliveries while the other sampled sets are shown in the Appendix (\Cref{tab:N78_bench}).

Our evaluation framework considers the following parameters:
\begin{enumerate}
    \item The number of available drones is fixed to $m = 10$;
    \item The number $N$ of deliveries is randomly chosen, between $10$ and $12$ for generating the dataset of larger problems, and between $7$ and $8$ for producing the smaller ones;
    \item The battery budget $B$ is randomly chosen between $50, 70$ and $100$;
    \item An $N$-dimensional list containing a random distribution of energy costs $c_j$ for completing each delivery $j \in N$, such that $c_j < B \ \forall \ j \in N$. To generate the lists, two configurations have been employed: a single Gaussian centered in $B/2$, and a uniform distribution between $0$ and $B$;
    \item An $N$-dimensional list of time windows $I_j$ within which to deliver item $j \in N$. Intervals are generated to extend from 8 a.m. to 8 p.m. and to be one or two hours long.
\end{enumerate}

In the subsequent subsection, comprehensive explanations are provided for the hardware utilized by each solver used on this benchmark.

\subsection{Annealers hardware}\label{sec:hardware}

By importing the necessary modules from the D-Wave library, it is possible to instantiate a solver object, both D-Wave's Simulated Annealer and the Quantum Annealer samplers. Then, one can submit the problem to the chosen sampler through the sample function, associated with the desired parameters, and retrieve results. We provide here the technical specifications of both the solver's hardware used. 

Simulated Annealing operates locally on personal computers without necessitating data transmission to specialized hardware. In particular, we used the Python environment along with D-Wave's Ocean SDK on our personal computer, equipped with an Intel(R) Core(TM) i7-11390H processor. Differently, performing experiments on Quantum Annealing requires specific hardware exploiting superconducting technology. Thus, to execute a quantum machine instruction (QMI), again through dedicated libraries within D-Wave's Ocean SDK, the information is sent across a network to the SAPI server, joining a queue for the chosen solver (\cite{documentation}). Our tests are conducted on the latest Quantum Annealer QPU hardware by D-Wave Advantage 1. The Advantage architecture counts 5627 qubits and uses the Pegasus graph topology, increasing the per-qubit connections to 15 (\cite{adv, chaindoc}). The solver's protocol and parameters characterize how the problem is run and allow control of the annealing process (e.g. chain strength, annealing time, and annealing schedule \cite{documentation}). In our tests, quantum annealing was implemented according to the standard forward annealing protocol, with a defaulted single-sample annealing time of 20.0$\mu s$. The chain strength parameter instead, was tuned from default to the smallest value able to increase the solver's performance for our test case, see \Cref{sec:chainstrength}. Also for the simulated annealing tests, we used the defaulted annealing protocol, corresponding to a linear interpolation of the $\beta$ parameter, within a defaulted range of values computed based on the problem biases  \cite{neal_documentation, neal_sample}. The number of reads was always set to 1000 for both classical and quantum tests.

Detailed explanations and comments are given in the following subsection for each metric used to evaluate the solver's performance.

\section{Numerical results}\label{sec:optimization-performance}

\begin{table*}[ht]
\centering
\renewcommand{\arraystretch}{0.9} 
\setlength{\tabcolsep}{4pt} 
\begin{tabularx}{\textwidth}{|l|l|*{3}{X|}*{3}{X|}*{2}{X|}|l|}
\hline
\textbf{Inst} & \textbf{Avg time} & \multicolumn{3}{c|}{\textbf{Solution}}                                                                                                   & \multicolumn{3}{c|}{\textbf{Drones}}                                                                                                     & \multicolumn{2}{c|}{\textbf{{[}battery, time, deliveries{]}}}                  & \textbf{Variables} \\ \hline
              &                   & \multicolumn{1}{c|}{\textbf{Avg}} & \multicolumn{1}{c|}{\cellcolor[HTML]{38FFF8}\textbf{Best}} & \textbf{ILP} & \multicolumn{1}{c|}{\textbf{Avg}} & \multicolumn{1}{c|}{\cellcolor[HTML]{38FFF8}\textbf{Best}} & \textbf{ILP} & \multicolumn{1}{c|}{\textbf{Avg}}      & \cellcolor[HTML]{38FFF8}\textbf{Best} & \textbf{}       \\ \hline
1             & 3.92              & \multicolumn{1}{c|}{87.0}         & \multicolumn{1}{c|}{\cellcolor[HTML]{38FFF8}86}            & 82              & \multicolumn{1}{c|}{8.6}          & \multicolumn{1}{c|}{\cellcolor[HTML]{38FFF8}8}             & 7               & \multicolumn{1}{c|}{{[}0.9 0.8 1. {]}} & \cellcolor[HTML]{3166FF}{[}1 1 1{]}   & 180             \\ \hline
2             & 3.51              & \multicolumn{1}{c|}{84.4}         & \multicolumn{1}{c|}{\cellcolor[HTML]{38FFF8}82}            & 74              & \multicolumn{1}{c|}{7.9}          & \multicolumn{1}{c|}{\cellcolor[HTML]{38FFF8}7}             & 5               & \multicolumn{1}{c|}{{[}0.9 0.5 1. {]}} & \cellcolor[HTML]{38FFF8}{[}0 0 1{]}   & 180             \\ \hline
3             & 3.87              & \multicolumn{1}{c|}{126.0}        & \multicolumn{1}{c|}{\cellcolor[HTML]{38FFF8}124}           & 122             & \multicolumn{1}{c|}{9}            & \multicolumn{1}{c|}{\cellcolor[HTML]{38FFF8}8}             & 7               & \multicolumn{1}{c|}{{[}0.7 0.4 1. {]}} & \cellcolor[HTML]{38FFF8}{[}1 0 1{]}   & 200             \\ \hline
4             & 3.43              & \multicolumn{1}{c|}{87.8}         & \multicolumn{1}{c|}{\cellcolor[HTML]{38FFF8}86}            & 80              & \multicolumn{1}{c|}{9}            & \multicolumn{1}{c|}{\cellcolor[HTML]{38FFF8}8}             & 6               & \multicolumn{1}{c|}{{[}1.  0.6 1. {]}} & \cellcolor[HTML]{38FFF8}{[}1 0 1{]}   & 180             \\ \hline
5             & 3.27              & \multicolumn{1}{c|}{85.4}         & \multicolumn{1}{c|}{\cellcolor[HTML]{38FFF8}84}            & 80              & \multicolumn{1}{c|}{7.9}          & \multicolumn{1}{c|}{\cellcolor[HTML]{38FFF8}7}             & 5               & \multicolumn{1}{c|}{{[}0.9 0.6 1. {]}} & \cellcolor[HTML]{38FFF8}{[}1 0 1{]}   & 180             \\ \hline
6             & 3.83              & \multicolumn{1}{c|}{115.8}        & \multicolumn{1}{c|}{\cellcolor[HTML]{38FFF8}110}           & 108             & \multicolumn{1}{c|}{7.1}          & \multicolumn{1}{c|}{\cellcolor[HTML]{38FFF8}6}             & 6               & \multicolumn{1}{c|}{{[}0.9 0.2 1. {]}} & \cellcolor[HTML]{38FFF8}{[}1 0 1{]}   & 190             \\ \hline
7             & 3.89              & \multicolumn{1}{c|}{126.6}        & \multicolumn{1}{c|}{\cellcolor[HTML]{3166FF}124}           & 124             & \multicolumn{1}{c|}{9.3}          & \multicolumn{1}{c|}{\cellcolor[HTML]{38FFF8}8}             & 8               & \multicolumn{1}{c|}{{[}0.9 0.4 1. {]}} & \cellcolor[HTML]{38FFF8}{[}1 0 1{]}   & 200             \\ \hline
8             & 3.65              & \multicolumn{1}{c|}{119.4}        & \multicolumn{1}{c|}{\cellcolor[HTML]{3166FF}118}           & 118             & \multicolumn{1}{c|}{8.5}          & \multicolumn{1}{c|}{\cellcolor[HTML]{38FFF8}8}             & 7               & \multicolumn{1}{c|}{{[}0.9 0.1 1. {]}} & \cellcolor[HTML]{38FFF8}{[}1 0 1{]}   & 190             \\ \hline
9             & 3.3               & \multicolumn{1}{c|}{88.1}         & \multicolumn{1}{c|}{\cellcolor[HTML]{38FFF8}77}            & 82              & \multicolumn{1}{c|}{8}            & \multicolumn{1}{c|}{\cellcolor[HTML]{38FFF8}7}             & 6               & \multicolumn{1}{c|}{{[}0.5 0.6 0.2{]}} & \cellcolor[HTML]{38FFF8}{[}0 1 0{]}   & 180             \\ \hline
10            & 4.69              & \multicolumn{1}{c|}{125.6}        & \multicolumn{1}{c|}{\cellcolor[HTML]{3166FF}122}           & 122             & \multicolumn{1}{c|}{8.7}          & \multicolumn{1}{c|}{\cellcolor[HTML]{38FFF8}7}             & 7               & \multicolumn{1}{c|}{{[}0.7 0.7 1. {]}} & \cellcolor[HTML]{38FFF8}{[}0 1 1{]}   & 200             \\ \hline
11            & 3.26              & \multicolumn{1}{c|}{87.4}         & \multicolumn{1}{c|}{\cellcolor[HTML]{38FFF8}84}            & 82              & \multicolumn{1}{c|}{8.7}          & \multicolumn{1}{c|}{\cellcolor[HTML]{38FFF8}7}             & 7               & \multicolumn{1}{c|}{{[}0.9 0.8 1. {]}} & \cellcolor[HTML]{3166FF}{[}1 1 1{]}   & 180             \\ \hline
12            & 3.97              & \multicolumn{1}{c|}{124.2}        & \multicolumn{1}{c|}{\cellcolor[HTML]{3166FF}118}           & 118             & \multicolumn{1}{c|}{8.7}          & \multicolumn{1}{c|}{\cellcolor[HTML]{38FFF8}7}             & 7               & \multicolumn{1}{c|}{{[}0.9 0.1 1. {]}} & \cellcolor[HTML]{38FFF8}{[}1 0 1{]}   & 200             \\ \hline
\end{tabularx}
\caption{Results over 10 attempts of SA on \eqref{eq:QUBO2} for the problem sets of \Cref{tab:k90_table} and with $k=120$ in \eqref{eq:constraint_all-at-once}. For each problem instance, we emphasize in light blue the "best solution", referring to the solution, among 10 results found by SA, that leads to the smallest value of $H_0$. When this value equals the exact solution, it is highlighted in dark blue. Moreover, the triplet of the constraint satisfaction is also colored in light blue, depicting 1 if the constraint is satisfied by the best solution, 0 otherwise. When each entry of the triplet is 1 (the best solution is feasible) it is colored in dark blue.}
\label{tab:findings}
\end{table*}

\begin{figure}[t!]
\begin{tikzpicture}
\begin{axis}[
    xlabel={$N$ (number of deliveries)},
    ylabel={Time (s)},
    ylabel near ticks,
    xlabel near ticks,
    xmin=4, xmax=8,
    ymin=0.15, ymax=1.6,
    xtick={4,5,6,7,8},
    ytick distance=0.1, 
    grid=both,
    axis lines=middle,
    enlargelimits=true,
    title={},
    legend style={at={(0.5,1.05)},anchor=south},
    ylabel style={xshift=-20pt, yshift=20pt},
    xlabel style={yshift=-30pt, xshift=-35pt} 
    ]

\addplot[
    only marks,
    mark=*,
    mark options={scale=0.8, fill=green}
] 
coordinates {
    (4,0.79)
    (5,0.86)
    (6,0.92)
    (7,1.05)
    (8,1.12)
};

\addplot[
    only marks,
    mark=*,
    mark options={scale=0.8, fill=blue}
] 
coordinates {
    (4,1.27)
    (5,1.42)
    (6,1.47)
    (7,1.64)
    (8,1.68)
};

\addplot[
    only marks,
    mark=*,
    mark options={scale=0.8, fill=yellow}
] 
coordinates {
    (4,0.136)
    (5,0.147)
    (6,0.154)
    (7,0.150)
    (8,0.156)
};

\addplot[
    only marks,
    mark=*,
    mark options={scale=0.8, fill=red}
] 
coordinates {
    (4,0.161)
    (5,0.168)
    (6,0.178)
    (7,0.164)
    (8,0.189)
};

\addlegendentry{SA on QUBO \eqref{eq:QUBO2}}
\addlegendentry{SA on QUBO \eqref{eq:QUBO1}}
\addlegendentry{QA on QUBO \eqref{eq:QUBO2}}
\addlegendentry{QA on QUBO \eqref{eq:QUBO1}}

\addplot[
    green,
    thick,
    mark=none
] 
coordinates {
    (4,0.79)
    (5,0.86)
    (6,0.92)
    (7,1.05)
    (8,1.12)
};

\addplot[
    blue,
    thick,
    mark=none
] 
coordinates {
    (4,1.27)
    (5,1.42)
    (6,1.47)
    (7,1.64)
    (8,1.68)
};

\addplot[
    yellow,
    thick,
    mark=none
] 
coordinates {
    (4,0.136)
    (5,0.147)
    (6,0.154)
    (7,0.150)
    (8,0.156)
};

\addplot[
    red,
    thick,
    mark=none
] 
coordinates {
    (4,0.161)
    (5,0.168)
    (6,0.178)
    (7,0.164)
    (8,0.189)
};

\end{axis}
\end{tikzpicture}
\caption{Time-to-solution comparison between SA and QA as a function of the number of deliveries $N$, when
applied to \cref{eq:QUBO1} and \cref{eq:QUBO2}. We let $N$ vary in $\{4,...,8\}$, we chose a fixed $m = 5$ and a fixed battery budget $B = 5$. Costs $c_j$ for each delivery $j \in N$ are randomly sampled from a Gaussian distribution such that $0 \le c_j \le
B \forall j \in N$. Time intervals $I_j$ are sampled from a uniform distribution from 8 a.m. to 8 p.m. and are at
least one hour long and at most two.}
\label{fig:time}
\end{figure}
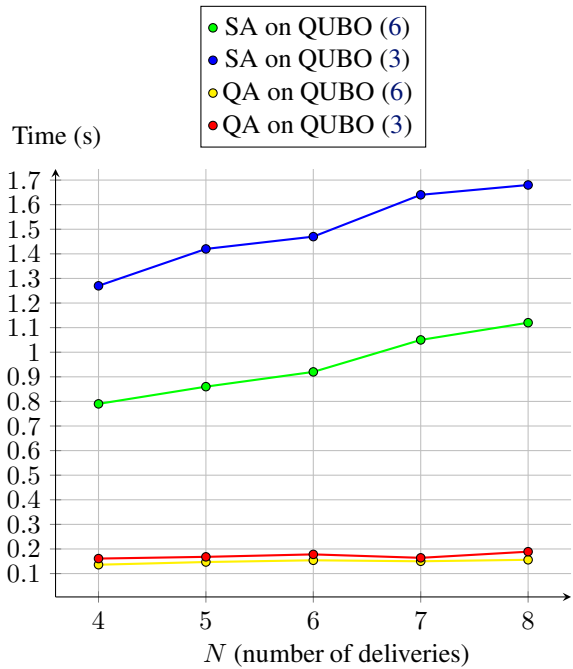

In this section, we provide a performance analysis for both quantum and simulated annealing samplers as compared to the results obtained using the original constrained ILP formulation \eqref{eq:ILP2}, solved with the commercial software Gurobi \cite{gurobi}. 

Being probabilistic methods, annealing techniques are employed iteratively, i.e.\ multiple runs of the optimization process are required to explore different regions of the solution space. 
Thus, for an annealing technique to be considered a successful probabilistic optimization method, it is sufficient to find at least one optimal and feasible solution among all the runs. Namely, the key metric for evaluating the performance of an annealing technique is not whether it consistently finds the optimal solution in every run, but rather if it can achieve the optimal solution at least once within multiple attempts. We seek the problem dimension range within which both Simulated and Quantum annealing results meet this requirement.

\subsection{Models comparison in time to solution}

The experiments detailed in the present section all pertain to the QUBO formulation in \Cref{eq:QUBO2}, as it was discussed in \Cref{qubomodels}
and \Cref{qaperf} to be significantly more convenient.

Also, we show in this subsection that the enhancement gained with \Cref{eq:QUBO2} over \Cref{eq:QUBO1} can be evaluated in terms of time-to-solution benefit. \Cref{fig:time} displays the impact of the variables' improvement in terms of time-to-solution, both for simulated annealing and quantum annealing. Precisely, the considered time for SA measures the extension of the classical algorithm execution, that retrieves the best sample among 1000 reads. In this sense, the plot shows how long a single call to SA takes to retrieve 1000 samples and subsequently extract the best. The reported values are averages computed within 10 of such calls to the solver. 

To make a fair comparison, we plot the QPU access time for QA, neglecting the data transmission latency. We notice that for both SA and QA there is a sensible time reduction when using the modified QUBO formulation \eqref{eq:QUBO2}. The figure also highlights the significant speed-up of QA with respect to SA. The effectiveness of the QUBO model \eqref{eq:QUBO2} in terms of solution quality will be thoroughly investigated in \Cref{sec:experiments}.

Note that \Cref{fig:time} shows just the QPU access time as a comparison with SA execution time. Namely, we chose to omit internet latency and preparation times, and also the time required by the embedding. However, while for the considered benchmark the embedding time makes a negligible contribution, we point out that, for complex embeddings, it can appear to be significant. The observed speed-up of quantum hardware is encouraging and aligned with their promising capacities in overcoming local minima related to purely quantum effects, as exposed in \Cref{qaperf}.

\subsection{Simulated Annealer results}\label{sec:SAresults}

\begin{table*}[ht]
\renewcommand{\arraystretch}{0.9} 
\setlength{\tabcolsep}{4pt} 
\begin{tabularx}{\textwidth}{|l|l|*{3}{X|}*{3}{X|}*{2}{X|}|l|}
\hline
\textbf{Inst} & \textbf{Avg time} & \multicolumn{3}{c|}{\textbf{Solution}}                                                                                                   & \multicolumn{3}{c|}{\textbf{Drones}}                                                                                                     & \multicolumn{2}{c|}{\textbf{{[}battery, time, deliveries{]}}}                  & \textbf{Variables} \\ \hline
              & \textbf{}         & \multicolumn{1}{c|}{\textbf{Avg}} & \multicolumn{1}{c|}{\cellcolor[HTML]{38FFF8}\textbf{Best}} & \textbf{ILP} & \multicolumn{1}{c|}{\textbf{Avg}} & \multicolumn{1}{c|}{\cellcolor[HTML]{38FFF8}\textbf{Best}} & \textbf{ILP} & \multicolumn{1}{c|}{\textbf{Avg}}      & \cellcolor[HTML]{38FFF8}\textbf{Best} &                    \\ \hline
1             & 2.64              & \multicolumn{1}{c|}{40.6}         & \multicolumn{1}{c|}{\cellcolor[HTML]{3166FF}38}            & 38              & \multicolumn{1}{c|}{6.3}          & \multicolumn{1}{c|}{\cellcolor[HTML]{38FFF8}5}             & 5               & \multicolumn{1}{c|}{{[}0.9 1.  1. {]}} & \cellcolor[HTML]{3166FF}{[}1 1 1{]}   & 150                \\ \hline
2             & 2.77              & \multicolumn{1}{c|}{54.8}         & \multicolumn{1}{c|}{\cellcolor[HTML]{38FFF8}52}            & 50              & \multicolumn{1}{c|}{7}            & \multicolumn{1}{c|}{\cellcolor[HTML]{38FFF8}6}             & 5               & \multicolumn{1}{c|}{{[}1.  0.9 0.8{]}} & \cellcolor[HTML]{3166FF}{[}1 1 1{]}   & 160                \\ \hline
3             & 7.58              & \multicolumn{1}{c|}{40.0}         & \multicolumn{1}{c|}{\cellcolor[HTML]{3166FF}40}            & 40              & \multicolumn{1}{c|}{6}            & \multicolumn{1}{c|}{\cellcolor[HTML]{38FFF8}6}             & 6               & \multicolumn{1}{c|}{{[}1. 0. 1.{]}}    & \cellcolor[HTML]{38FFF8}{[}1 0 1{]}   & 140                \\ \hline
4             & 5.48              & \multicolumn{1}{c|}{41.4}         & \multicolumn{1}{c|}{\cellcolor[HTML]{38FFF8}38}            & 36              & \multicolumn{1}{c|}{5.9}          & \multicolumn{1}{c|}{\cellcolor[HTML]{38FFF8}5}             & 4               & \multicolumn{1}{c|}{{[}0.9 1.  0.6{]}} & \cellcolor[HTML]{3166FF}{[}1 1 1{]}   & 150                \\ \hline
5             & 2.72              & \multicolumn{1}{c|}{54.0}         & \multicolumn{1}{c|}{\cellcolor[HTML]{3166FF}54}            & 54              & \multicolumn{1}{c|}{7}            & \multicolumn{1}{c|}{\cellcolor[HTML]{38FFF8}7}             & 7               & \multicolumn{1}{c|}{{[}1. 1. 1.{]}}    & \cellcolor[HTML]{3166FF}{[}1 1 1{]}   & 150                \\ \hline
6             & 2.73              & \multicolumn{1}{c|}{39.6}         & \multicolumn{1}{c|}{\cellcolor[HTML]{38FFF8}34}            & 32              & \multicolumn{1}{c|}{6.1}          & \multicolumn{1}{c|}{\cellcolor[HTML]{38FFF8}4}             & 3               & \multicolumn{1}{c|}{{[}0.9 0.9 1. {]}} & \cellcolor[HTML]{38FFF8}{[}0 1 1{]}   & 150                \\ \hline
7             & 2.48              & \multicolumn{1}{c|}{39.0}         & \multicolumn{1}{c|}{\cellcolor[HTML]{38FFF8}38}            & 34              & \multicolumn{1}{c|}{5.5}          & \multicolumn{1}{c|}{\cellcolor[HTML]{38FFF8}5}             & 4               & \multicolumn{1}{c|}{{[}0.8 1.  1. {]}} & \cellcolor[HTML]{38FFF8}{[}0 1 1{]}   & 140                \\ \hline
8             & 2.83              & \multicolumn{1}{c|}{54.9}         & \multicolumn{1}{c|}{\cellcolor[HTML]{3166FF}50}            & 50              & \multicolumn{1}{c|}{6.5}          & \multicolumn{1}{c|}{\cellcolor[HTML]{38FFF8}5}             & 5               & \multicolumn{1}{c|}{{[}0.5 0.7 0.7{]}} & \cellcolor[HTML]{38FFF8}{[}0 0 1{]}   & 160                \\ \hline
9             & 2.55              & \multicolumn{1}{c|}{40.4}         & \multicolumn{1}{c|}{\cellcolor[HTML]{3166FF}38}            & 38              & \multicolumn{1}{c|}{6}            & \multicolumn{1}{c|}{\cellcolor[HTML]{38FFF8}5}             & 5               & \multicolumn{1}{c|}{{[}1.  0.9 0.9{]}} & \cellcolor[HTML]{3166FF}{[}1 1 1{]}   & 150                \\ \hline
10            & 2.83              & \multicolumn{1}{c|}{54.0}         & \multicolumn{1}{c|}{\cellcolor[HTML]{38FFF8}52}            & 50              & \multicolumn{1}{c|}{7}            & \multicolumn{1}{c|}{\cellcolor[HTML]{38FFF8}6}             & 5               & \multicolumn{1}{c|}{{[}0.9 0.9 1. {]}} & \cellcolor[HTML]{3166FF}{[}1 1 1{]}   & 160                \\ \hline
11            & 2.7               & \multicolumn{1}{c|}{39.4}         & \multicolumn{1}{c|}{\cellcolor[HTML]{38FFF8}38}            & 32              & \multicolumn{1}{c|}{5.7}          & \multicolumn{1}{c|}{\cellcolor[HTML]{38FFF8}5}             & 3               & \multicolumn{1}{c|}{{[}1.  0.7 1. {]}} & \cellcolor[HTML]{38FFF8}{[}1 0 1{]}   & 150                \\ \hline
12            & 3.6               & \multicolumn{1}{c|}{36.8}         & \multicolumn{1}{c|}{\cellcolor[HTML]{3166FF}36}            & 36              & \multicolumn{1}{c|}{5.1}          & \multicolumn{1}{c|}{\cellcolor[HTML]{38FFF8}5}             & 5               & \multicolumn{1}{c|}{{[}1.  0.2 1. {]}} & \cellcolor[HTML]{38FFF8}{[}1 0 1{]}   & 140                \\ \hline
\end{tabularx}

\caption{Results over 10 attempts of SA on \eqref{eq:QUBO2} for the problem sets of \Cref{tab:N78_bench} and with $k=120$ in \eqref{eq:constraint_all-at-once}. For each problem instance, we emphasize in light blue the "best solution", referring to the solution, among 10 results found by SA, that leads to the smallest value of $H_0$. When this value equals the exact solution, it is highlighted in dark blue. Moreover, the triplet of the constraint satisfaction is also colored in light blue, depicting 1 if the constraint is satisfied by the best solution, and 0 otherwise. When each entry of the triplet is 1 (the best solution is feasible) it is colored in dark blue.}
\label{tab:findings2}
\end{table*}

The results found by SA on each of the 24 synthetically generated problem instances reported in \Cref{tab:k90_table} and \Cref{tab:N78_bench}, are shown in \Cref{tab:findings} and \Cref{tab:findings2} respectively. There, we show both average performance and "best" performance, obtained by showing the results corresponding to the smallest value of $H_0$, over 10 runs with starting points randomly chosen. As discussed in \Cref{sec:optimization-performance}, the results for the "best" run serve as a crucial indicator of whether the solver has successfully located the global minimum or whether it has gotten stuck into a local minimum, even after multiple attempts.

Overall, \Cref{tab:findings} and \Cref{tab:findings2} summarize the results obtained across all the 24 problem instances described in \Cref{sec:problems-setup}, providing details about the following metrics:

\begin{itemize}
    \item Avg time: the time to solution required by the SA to retrieve the best sample for each of the 12 instances. This value is an average computed within 10 calls to the solver. Precisely, the considered time for SA measures the extension of the classical algorithm execution, that retrieves the best sample among 1000 reads;
    \item Solution, avg: the average value, within 10 calls to the SA, of the best sample retrieved by the solver. In each call, after the sampling process, the solver extracts the best value of the objective function among the various solutions found; 
    \item Solution, best: in agreement with what has been said, the key information is whether the solver finds the real minimum at least once. In order to understand this, for each of the 12 instances, we select the smallest among the 10 results found by SA. We highlight in darker blue the cases in which it corresponds to the real global minimum;
    \item Solution, ILP: the correct optimal value of the objective function for each instance, computed with a commercial deterministic solver branch-and-bound based \cite{gurobi}, applied to the original constrained ILP \eqref{eq:ILP1}. 
    \item Drones, avg: the number of drones corresponding to the solution with the smallest value of the objective function, for each instance. We report here the average number of such drones over 10 calls to SA;
    \item Drones, best: for each of the 12 instances, the number of drones corresponding to the solution with the smallest objective function value. We highlight in darker blue the cases in which it corresponds to the real minimum number of prescribed drones;
    \item Drones, ILP: the number of drones prescribed by the exact solution given by the deterministic solver;
    \item $[$battery, time-consistency, all deliveries once$]$, avg: triplet of percentages of satisfaction, namely how many times, in 10 attempts, each of the specific constraints is satisfied; 
    \item $[$battery, time-consistency, all deliveries once$]$, \\ best: triplet of satisfaction, specific for the solution with the smallest objective function value. The triplet entry equals 1 if the corresponding constraint is satisfied, and 0 otherwise.   
\end{itemize}

The results highlighted with a dark blue color in \Cref{tab:findings} and \Cref{tab:findings2} emphasize the capacity of the solver to find both optimal and feasible solutions at least once. We recall that by "feasible solution" we mean a solution that satisfies simultaneously all the problem constraints, namely bringing to completion the required delivery task. It is also important to stress that even though the simultaneous satisfaction of the constraints is an important feature, studying how the method addresses single constraints is also useful. In fact, it allows us to understand potential disparities among different constraints, thus hinting at the necessity of tuning the corresponding penalty weight, as explained in \Cref{penalties}.

\Cref{tab:findings} shows that utilizing a benchmark comprising 10 or 12 deliveries, where each instance consists of approximately 200 variables, SA struggles to find solutions that are both feasible and optimal.
However, the situation changes significantly when lowering the problem dimension. When focusing on 7 or 8 deliveries, encompassing approximately 150 variables, \Cref{tab:findings2} shows that SA is able to find solutions that are both feasible and optimal. The comparison of \Cref{tab:findings} and \Cref{tab:findings2} highlights the dependence of SA performance on the number of variables. It underscores that variable reduction can play a key role, enabling the method to find satisfactory solutions. More specifically, we identify a dimension-limiting case of $\sim 150$ variables within which, with this setting and these parameters, the solver can be considered comparable to deterministic methods.

In the following subsection, we aim to conduct a similar analysis to the one undertaken for SA, but for the QPU described in \Cref{sec:hardware}. However, we will see that for the quantum case, due to hardware constraints, the dimension sensitivity is impacting the annealing performance even more. 

\subsection{Quantum Annealer results}\label{qa}

\begin{table}[htbp]
\begin{tabularx}{\columnwidth}{|c|X|c|c|}
\hline
\textbf{N} & \begin{small}\textbf{[battery,time,dels]}\end{small} & \textbf{\% Feasibility} & \textbf{Qubits} \\ \hline
4          & [0.8 1. 0.7]                         & 0.7                     & 109             \\ \hline
5          & [0.8 0.8 1.]                         & 0.7                     & 115             \\ \hline
6          & [0.8 0.2 0.6]                        & 0.2                     & 158             \\ \hline
7          & [0.6 0.4 0.2]                        & 0.0                     & 193             \\ \hline
8          & [0.6 0.9 0.1]                        & 0.0                     & 222             \\ \hline
\end{tabularx}
\caption{Constraints satisfaction and feasibility check with respect to number $N$ of deliveries and relative qubits, with $m = 4$ and $k = 120$.}
\label{tab:feasibility}
\end{table}

We address the issues related to the problem embedding into the quantum hardware described in \Cref{emb}, for the specific problem at hand, and how these influence the performance sensitivity concerning problem size and chain strength. The goal is to identify, as done with the classical annealer, a dimension range of applicability of currently available machines for the DDPP formulated as in \Cref{eq:QUBO2}.

Using QA for the bigger benchmark problems in \Cref{tab:k90_table} requires $\sim 800$ physical qubits, yielding poorly performing results. In general, encountering many infeasible solutions represents a notable drawback for annealing techniques. This underscores the need for establishing a dimensionality range within which these methods, given current technological capabilities, appear to be reliable. \Cref{tab:findings2} demonstrates the capacity of SA to derive meaningful solutions for the DDPP on instances featuring $\sim 150$ variables. Consequently, we seek a similar estimation for QA: we aim to identify the largest problem size for which QA is able to find at least one feasible solution while approximating the optimal value of $H_0$.

Starting from \Cref{tab:k90_table} data, we repeat the same experiments done with SA, changing the solver to QA, and significantly decreasing $m$ and $N$. All tests are done according to \Cref{sec:SAresults} setup, just changing the solver and the problem size. Particularly, 10 QA runs with random starting points are done for each decreased-size instance, each with 1000 reads. For what concerns the size reduction, we consider a fixed $m = 4$ and tune $N$ from 4 to the largest value that results in at least one feasible solution, determined to be $N=6$. The detailed findings are presented in \Cref{tab:feasibility}, revealing that on this simpler use-case, QA can be considered successful.

As the number of deliveries increased, so did the likelihood of encountering unfeasible solutions. Larger problem instances amplify the complexity of exploring the solution space, particularly due to the presence of multiple constraints. Our tests show that QA is robust up to a threshold of approximately 160 physical variables. This corresponds to scheduling $m=4$ drones on bringing to completion $N=6$ deliveries, with prescribed constraints of time and battery. Beyond this point, the performance of QA declines. The next subsection explores the possibility of stretching this limiting problem dimension through parameter tuning.

\subsection{Chain Strength Tuning}\label{sec:chainstrength}

\begin{table}[ht]
\small
\begin{tabularx}{\columnwidth}{|c|c|X|}
\hline
\multicolumn{3}{|c|}{m = 4, N = 7, k = 120}                                                                                                               \\ \hline
\textbf{Chain Strength}                  & \textbf{{[}battery, time, dels{]}}  & \textbf{\%Feas} \\ \hline
Default + 1 000                          & [0.4 0.4 0.1]                         & 0.0                     \\ \hline
\rowcolor[HTML]{CBCEFB} 
\cellcolor[HTML]{CBCEFB}Default + 5 000  & \cellcolor[HTML]{CBCEFB}[0.4 0.5 0.2] & 0.1                     \\ \hline
Default + 8 000                          & [0.4 0.4 0.3]                         & 0.0                     \\ \hline
\rowcolor[HTML]{9698ED} 
\cellcolor[HTML]{9698ED}Default + 10 000 & \cellcolor[HTML]{9698ED}[0.7 0.6 0.4] & 0.2                     \\ \hline
\cellcolor[HTML]{9698ED}Default + 20 000 & \cellcolor[HTML]{9698ED}[0.6 0.5 0.6] & \cellcolor[HTML]{9698ED} 0.2                     \\ \hline
\cellcolor[HTML]{9698ED}Default + 30 000 & \cellcolor[HTML]{9698ED}[0.6 0.4 0.4] & \cellcolor[HTML]{9698ED} 0.2                     \\ \hline
\cellcolor[HTML]{6665CD}Default + 40 000 & \cellcolor[HTML]{6665CD}[0.9 0.6 0.5] & \cellcolor[HTML]{6665CD} 0.3                     \\ \hline
\rowcolor[HTML]{CBCEFB} 
\cellcolor[HTML]{CBCEFB}Default + 50 000 & \cellcolor[HTML]{CBCEFB}[0.6 0.6 0.4] & 0.1                     \\ \hline
Default + 100 000                        & [0.3 0.2 0.4]                         & 0.0                     \\ \hline
\end{tabularx}
\caption{Chain Strength tuning
for threshold-size case of $N = 7, m = 4$. Default chain strength value
for the indicated case corresponds to 109413.47.}
\label{tab:cs}
\end{table}

In this subsection, we observe that the dimension threshold can be raised from $N=6$ to $N=7$ tuning the chain strength of the annealer. To this end, we look for the optimal choice of this parameter such that QA yields at least one feasible solution to the smallest non-functioning case: $m=4, N=7$. Larger choices of the parameters turned out to be untreatable by the available quantum hardware. 

We concentrate on the chain-strength parameter \cite{chainstrength} arising specifically when utilizing QA technology. This value tunes the coupling force applied to chain qubits inside the hardware, in order to faithfully represent single QUBO variables, as explained in \Cref{emb}. In other words, couplers are assigned a bias that favors chain qubits taking the same value, and this is controlled by the chain strength parameter.

Knowing the threshold sizes of $N = 7, m = 4$ for the QA to approach the DDPP successfully, we manually increase the chain strength parameter corresponding to this indicted case, in a tuning process. Indeed, given the default chain strength value of $109413.47$, for this size, we iteratively raise it, and report in \Cref{tab:cs} how this impacts solution quality in terms of feasibility. Parameters that lead to a feasibility ratio boost are highlighted in  \Cref{tab:cs}, announcing the possibility of raising the reliability of QA on this problem up to $N=7$ deliveries. This shows that this case, which corresponds to a large number of physical qubits (between 190 and 200), can be addressed by currently available quantum annealers by optimizing the chain strength.

Finally, to provide a more complete overview of the performance of QA, in the next subsection, we tackle the issue of scalability and show how the performance changes when increasing $N$ in the derived range.

\subsection{Performance scaling with respect to the number of deliveries}

In this section, we address the scalability of QA and SA solvers, by examining how performance varies within the approachable size of $N \in \{4,5,6,7\}$. Specifically, a study of the sensitivity of the annealers' solutions with respect to the number of deliveries is given and evaluated through a comparison with deterministic solver findings on the same problems.

We aim to tune $N \in \{4,5,6,7,8\}$ since, as we recall from \Cref{qa}, QA performs successfully on all values of $N \in \{4,5,6,7\}$, when the set of available drones has dimension $m=4$. Therefore, for conducting the $N$-tuning in the stated range, for QA tests, we fix $m=4$. We chose to show also the case with $N=8$ deliveries, as representative of the dimension-limiting case for QA. Similarly, SA handles problems in the stated $N$ range for $m=10$, as shown in \Cref{sec:optimization-performance}. Therefore, for the same tests with SA, we fix $m=10$. For each $N$ value and the fixed $m$ values, we generate a test problem, as done for the instances from \Cref{tab:k90_table}. Of course, these synthetic instances, created for each $N$ value, comprise $N$-dimensional lists of energy costs and time intervals, and consequently are of increasing size, as $N$ increases. Specific sampled sets of parameters are shown in \Cref{tab:incr_N_bench} of the Appendix.

We report in \Cref{fig:qavssa} the results for SA and QA on these random instances of increasing size. 
In the figure, we report the average optimal objective function values over the feasible solutions obtained with 10 calls to the annealers. 
These values are compared with the corresponding values obtained by the commercial deterministic solver applied to the ILP formulation  \eqref{eq:ILP2}. We also show the percentage of feasible solutions found in 10 attempts, for each value of $N$ and for both SA and QA. 
As we average only over feasible solutions, these figures provide information on both the optimality of the computed output (as compared to the deterministic solution) and its feasibility.

According to our findings, with the penalty weights specified in \Cref{penalties} and the optimal chain strength for the $N=7$ case of \Cref{tab:cs}, QA is effective on the selected size range, finding meaningful solutions. In terms of physical qubits, this shows that QA can competitively tackle DDPP instances, as long as the number of variables remains under $\sim$ 190. This is a noteworthy value and is obtained due to the chain strength tuning in \Cref{tab:cs}. 

\begin{figure}[ht]
\begin{tikzpicture}

\begin{axis}[
    xlabel={$N$},
    ylabel={Energy},
    xmin=4, xmax=8,
    ymin=5, ymax=55,
    xtick={4,5,6,7,8},
    ytick distance=5, 
    grid=both,
    axis lines=middle,
    enlargelimits=true,
    title={QA on QUBO \eqref{eq:QUBO2}, $m = 4$},
    legend style={at={(0.15,0.65)},anchor=south}
]

\addplot[
    only marks,
    mark=x,
    mark options={scale=2.5, fill=black}
] 
coordinates {
    (4,8)
    (5,16)
    (6,26)
    (7,36)
    (8,50)
};

\addplot[
    only marks,
    mark=*,
    mark options={scale=1.5, fill=red}
] 
coordinates {
    (4,10.8)
    (5, 18)
    (6,26)
    (7,36)
};

\addlegendentry{Exact}
\addlegendentry{QA}

\addplot[
    red,
    thick,
    mark=none
] 
coordinates {
    (4,10.8)
    (5, 18)
    (6,26)
    (7,36)
};

\addplot[
    blue,
    thick,
    mark=none
] 
coordinates {
    (4,10)
};

\addplot[
    black,
    thick,
    mark=none
] 
coordinates {
    (4,8)
    (5,16)
    (6,26)
    (7,36)
    (8,50)
};

\node[red, anchor=south] at (axis cs:4,12) {70\%};
\node[red, anchor=south] at (axis cs:5,20) {70\%};
\node[red, anchor=south] at (axis cs:6,28) {20\%};
\node[red, anchor=south] at (axis cs:7,38) {30\%};
\end{axis}
\end{tikzpicture}

\begin{tikzpicture}
\begin{axis}[
    xlabel={$N$},
    ylabel={Energy},
    xmin=4, xmax=8,
    ymin=5, ymax=55,
    xtick={4,5,6,7,8},
    ytick distance=5, 
    grid=both,
    axis lines=middle,
    enlargelimits=true,
    title={SA on QUBO \eqref{eq:QUBO2}, $m = 10$},
    legend style={at={(0.15,0.65)},anchor=south}
]
\addplot[
    only marks,
    mark=x,
    mark options={scale=2.5, fill=black}
] 
coordinates {
    (4,10)
    (5,16)
    (6,30)
    (7,38)
    (8,50)
};

\addplot[
    only marks,
    mark=*,
    mark options={scale=1.5, fill=blue}
] 
coordinates {
    (4,11.8)
    (5,18)
    (6,29.6)
    (7,39.5)
    (8,53.3)
};

\addlegendentry{Exact}
\addlegendentry{SA}

\addplot[
    blue,
    thick,
    mark=none
] 
coordinates {
    (4,11.8)
    (5,18)
    (6,29.6)
    (7,39.5)
    (8,53.3)
};

\addplot[
    black,
    thick,
    mark=none
] 
coordinates {
    (4,10)
    (5,16)
    (6,30)
    (7,38)
    (8,50)
};

\node[blue, anchor=south] at (axis cs:4,13.8) {100\%};
\node[blue, anchor=south] at (axis cs:5,20) {80\%};
\node[blue, anchor=south] at (axis cs:6,31.6) {60\%};
\node[blue, anchor=south] at (axis cs:7,41.5) {90\%};
\node[blue, anchor=south] at (axis cs:8,55.3) {30\%};

\end{axis}

\end{tikzpicture}
\caption{Percentages of feasibility and average (over feasible solutions) of optimal values of the objective function found by applying QA and SA to the QUBO formulation \eqref{eq:QUBO2}. \texttt{Exact} corresponds to the solutions obtained with a deterministic solver applied to the ILP formulation \eqref{eq:ILP2}.}
\label{fig:qavssa}
\end{figure}
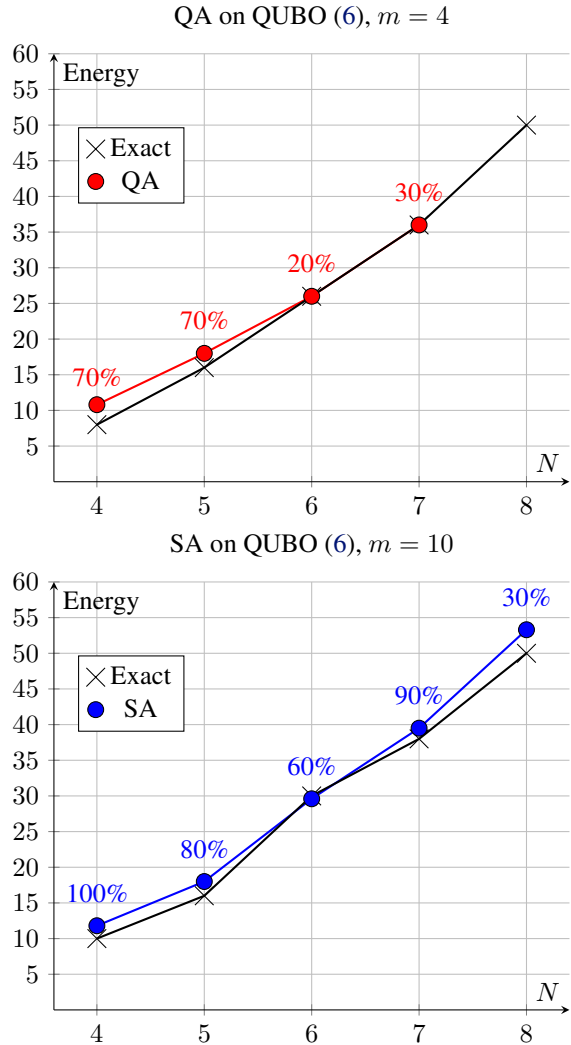

\section{Conclusion and future work}

The overarching goal of the work was to study the extent to which annealing techniques can be used to tackle the DDPP, a popular and notoriously complex optimization problem.

We presented the standard QUBO formulation of the problem and proposed a novel one leading to a significant reduction in constraints and thus in the variable count. We conducted a series of experiments with the primary objective of exploring the possibility of using quantum annealing machines and understanding the correct dimension range within which QA works well. We undertook an analysis involving both the Simulated and the Quantum Annealer solvers by D-Wave.

Our study of the potential applicability of quantum annealing hardware is conceived with a forward-looking perspective, as technology advances and yields larger and more interconnected devices. The analysis of how QA scales using small problems is particularly relevant also in view of the potential application within hybrid quantum-classical computing frameworks, where the quantum component tackles specific sub-problems of larger classical algorithms \cite{hybrid2, optimizingtheoptimizer, hybrid}. It is a topic of interest to develop strategies for handling larger problems with QA alone, as well as in cooperation with other methods.

We observe from our study that, even though QA can efficiently approximate exact solutions, still commercial deterministic solvers show the best performance. However, it is also known that, for big-size instances, exact solvers like Gurobi become unable to obtain optimal solutions with zero optimality gap, within a reasonable amount of time. For example, in \cite{hybrid2} it is shown that the branch-and-cut algorithm implemented in Gurobi was unable to solve the considered optimization problem within 24 hours, and hence yielded a sub-optimal solution. In that work, a hybrid quantum-classical partitioning approach instead was shown to find good quality near-optimal solutions within reasonable computational time. Therefore, depending on the problem setting, QA can represent a valid alternative in a hybrid-decomposition framework, as opposed to approaching the global problem with a deterministic solver. 

\section*{Acknowledgments}
We acknowledge the CINECA award under the ISCRA initiative, for the availability of quantum computing resources and support. S.T.\ acknowledges support from Leonardo SpA and wishes to thank the Quantum and HPC Labs in Genova for their hospitality.  

\bibliography{references}

\begin{thebibliography}{37}
\providecommand{\natexlab}[1]{#1}
\providecommand{\url}[1]{\texttt{#1}}
\expandafter\ifx\csname urlstyle\endcsname\relax
  \providecommand{\doi}[1]{doi: #1}\else
  \providecommand{\doi}{doi: \begingroup \urlstyle{rm}\Url}\fi

\bibitem[Ajagekar et~al.(2020)Ajagekar, Humble, and You]{hybrid2}
Ajagekar, A., Humble, T., and You, F.
\newblock Quantum computing based hybrid solution strategies for large-scale
  discrete-continuous optimization problems.
\newblock \emph{Computers amp; Chemical Engineering}, 132:\penalty0 106630,
  January 2020.
\newblock ISSN 0098-1354.
\newblock \doi{10.1016/j.compchemeng.2019.106630}.
\newblock URL \url{http://dx.doi.org/10.1016/j.compchemeng.2019.106630}.

\bibitem[Ayodele(2022)]{penalties}
Ayodele, M.
\newblock \emph{Penalty Weights in QUBO Formulations: Permutation Problems},
  pp.\  159–174.
\newblock Springer International Publishing, 2022.
\newblock ISBN 9783031041488.
\newblock \doi{10.1007/978-3-031-04148-8_11}.
\newblock URL \url{http://dx.doi.org/10.1007/978-3-031-04148-8_11}.

\bibitem[Bass et~al.(2020)Bass, Henderson, Heath, and
  au2]{optimizingtheoptimizer}
Bass, G., Henderson, M., Heath, J., and au2, J. D.~I.
\newblock Optimizing the optimizer: Decomposition techniques for quantum
  annealing, 2020.

\bibitem[Bergermann et~al.(2024)Bergermann, Stoll, and
  Tudisco]{bergermann2024nonlinear}
Bergermann, K., Stoll, M., and Tudisco, F.
\newblock A nonlinear spectral core-periphery detection method for multiplex
  networks.
\newblock \emph{Proceedings Royal Society A}, 2024.

\bibitem[Bettina~Heim \& Troyer(2015)Bettina~Heim and
  Troyer]{qaoutperformance2}
Bettina~Heim, Troels F.~Rønnow, S. V.~I. and Troyer, M.
\newblock Quantum versus classical annealing of ising spin glasses.
\newblock \emph{Science}, 2015.
\newblock \doi{10.1126/science.aaa4170}.

\bibitem[Boixo et~al.(2018)Boixo, Isakov, Smelyanskiy, Babbush, Ding, Jiang,
  Bremner, Martinis, and Neven]{Boixo}
Boixo, S., Isakov, S.~V., Smelyanskiy, V.~N., Babbush, R., Ding, N., Jiang, Z.,
  Bremner, M.~J., Martinis, J.~M., and Neven, H.
\newblock Characterizing quantum supremacy in near-term devices.
\newblock \emph{Nature Physics}, 14\penalty0 (6), 2018.
\newblock ISSN 1745-2481.
\newblock \doi{10.1038/s41567-018-0124-x}.
\newblock URL \url{http://dx.doi.org/10.1038/s41567-018-0124-x}.

\bibitem[Bozhedarov et~al.(2024)Bozhedarov, Usmanov, Salakhov, Boev, Kiktenko,
  and Fedorov]{bppfuel}
Bozhedarov, A.~A., Usmanov, S.~R., Salakhov, G.~V., Boev, A.~S., Kiktenko,
  E.~O., and Fedorov, A.~K.
\newblock Quantum and quantum-inspired optimization for solving the minimum bin
  packing problem.
\newblock \emph{Journal of Physics: Conference Series}, 2701\penalty0
  (1):\penalty0 012129, February 2024.
\newblock ISSN 1742-6596.
\newblock \doi{10.1088/1742-6596/2701/1/012129}.
\newblock URL \url{http://dx.doi.org/10.1088/1742-6596/2701/1/012129}.

\bibitem[Cai et~al.(2014)Cai, Macready, and Roy]{minors}
Cai, J., Macready, W.~G., and Roy, A.
\newblock A practical heuristic for finding graph minors, 2014.

\bibitem[Crosson \& Lidar(2021)Crosson and Lidar]{qaenhancement}
Crosson, E.~J. and Lidar, D.~A.
\newblock Prospects for quantum enhancement with diabatic quantum annealing.
\newblock \emph{Nature Reviews Physics}, 2021.
\newblock \doi{10.1038/s42254-021-00313-6}.

\bibitem[{D-Wave Systems Inc.}(2020{\natexlab{a}})]{chaindoc}
{D-Wave Systems Inc.}
\newblock D-wave qpu architecture: Topologies.
\newblock Webpage, 2020{\natexlab{a}}.
\newblock URL
  \url{https://docs.dwavesys.com/docs/latest/c_gs_4.html#pegasusinternalcoupled}.
\newblock Accessed on 2024-02-29.

\bibitem[{D-Wave Systems Inc.}(2020{\natexlab{b}})]{chainstrength}
{D-Wave Systems Inc.}
\newblock Programming the d-wave qpu: Setting the chain strength.
\newblock Webpage, 2020{\natexlab{b}}.
\newblock URL
  \url{https://www.dwavesys.com/media/vsufwv1d/14-1041a-a_setting_the_chain_strength.pdf}.
\newblock Accessed on 2024-03-08.

\bibitem[{D-Wave Systems Inc.}(2021)]{documentation}
{D-Wave Systems Inc.}
\newblock D-wave-system documentation release 1.10.0.
\newblock Webpage, 2021.
\newblock URL
  \url{https://docs.ocean.dwavesys.com/_/downloads/system/en/latest/pdf/}.
\newblock Accessed on 2024-03-08.

\bibitem[{D-Wave Systems Inc.}(2024)]{neal_sample}
{D-Wave Systems Inc.}
\newblock Simulatedannealingsampler.sample method, 2024.
\newblock URL
  \url{https://docs.ocean.dwavesys.com/en/stable/docs_neal/reference/generated/neal.sampler.SimulatedAnnealingSampler.sample.html}.
\newblock Accessed: 2024-06-05.

\bibitem[de~Andoin et~al.(2022)de~Andoin, Osaba, Oregi, Villar-Rodriguez, and
  Sanz]{bpp}
de~Andoin, M.~G., Osaba, E., Oregi, I., Villar-Rodriguez, E., and Sanz, M.
\newblock Hybrid quantum-classical heuristic for the bin packing problem.
\newblock In \emph{Proceedings of the Genetic and Evolutionary Computation
  Conference Companion}, GECCO ’22. ACM, July 2022.
\newblock \doi{10.1145/3520304.3533986}.
\newblock URL \url{http://dx.doi.org/10.1145/3520304.3533986}.

\bibitem[et. al.(2019)]{qaenhancement2}
et. al., F.~A.
\newblock Quantum supremacy using a programmable superconducting processor.
\newblock \emph{Nature}, 2019.
\newblock \doi{10.1038/s41586-019-1666-5}.

\bibitem[et. al.(2018)]{qaoutperformance4}
et. al., R.~H.
\newblock Phase transitions in a programmable quantum spin glass simulator.
\newblock \emph{Science}, 2018.
\newblock \doi{10.1126/science.aat2025}.

\bibitem[Feynman et~al.(1965)Feynman, Leighton, and Sands]{physics}
Feynman, R.~P., Leighton, R.~B., and Sands, M.
\newblock \emph{The Feynman lectures on physics / Vol. 3, Quantum mechanics}.
\newblock Addison-Wesley, Reading, Mass., 1965.
\newblock \doi{10.22331/idonotexist}.

\bibitem[Glover et~al.(2019)Glover, Kochenberger, and Du]{qubo1}
Glover, F., Kochenberger, G., and Du, Y.
\newblock A tutorial on formulating and using qubo models, 2019.

\bibitem[{Gurobi Optimizer}(2020)]{gurobi}
{Gurobi Optimizer}.
\newblock Gurobi documentation.
\newblock Online, 2020.
\newblock URL \url{https://www.gurobi.com/documentation/}.
\newblock 2024-02-29.

\bibitem[Harikrishnakumar et~al.(2020)Harikrishnakumar, Nannapaneni, Nguyen,
  Steck, and Behrman]{veh2}
Harikrishnakumar, R., Nannapaneni, S., Nguyen, N.~H., Steck, J.~E., and
  Behrman, E.~C.
\newblock A quantum annealing approach for dynamic multi-depot capacitated
  vehicle routing problem, 2020.

\bibitem[Higham et~al.(2022)Higham, Higham, and Tudisco]{higham2022testing}
Higham, C.~F., Higham, D.~J., and Tudisco, F.
\newblock Core-periphery partitioning and quantum annealing.
\newblock In \emph{Proceedings of the 28th ACM SIGKDD Conference on Knowledge
  Discovery and Data Mining}, pp.\  565--573, 2022.

\bibitem[Inc.(2024)]{neal_documentation}
Inc., D.-W.~S.
\newblock Neal documentation.
\newblock \url{https://docs.ocean.dwavesys.com/projects/neal/en/latest/}, 2024.
\newblock Accessed: 2024-06-06.

\bibitem[Irie et~al.(2019)Irie, Wongpaisarnsin, Terabe, Miki, and
  Taguchi]{veh4}
Irie, H., Wongpaisarnsin, G., Terabe, M., Miki, A., and Taguchi, S.
\newblock Quantum annealing of vehicle routing problem with time, state and
  capacity, 2019.

\bibitem[Jana \& Mandal(2023)Jana and Mandal]{bppjanamandal}
Jana, S. and Mandal, P.~S.
\newblock Approximation algorithms for drone delivery packing problem.
\newblock In \emph{Proceedings of the 24th International Conference on
  Distributed Computing and Networking}, ICDCN '23, pp.\  262–269, New York,
  NY, USA, 2023. Association for Computing Machinery.
\newblock ISBN 9781450397964.
\newblock \doi{10.1145/3571306.3571411}.
\newblock URL \url{https://doi.org/10.1145/3571306.3571411}.

\bibitem[Kadowaki \& Nishimori(1998)Kadowaki and Nishimori]{Kadowaki_1998}
Kadowaki, T. and Nishimori, H.
\newblock Quantum annealing in the transverse ising model.
\newblock \emph{Physical Review E}, 58\penalty0 (5):\penalty0 5355–5363,
  November 1998.
\newblock ISSN 1095-3787.
\newblock \doi{10.1103/physreve.58.5355}.
\newblock URL \url{http://dx.doi.org/10.1103/PhysRevE.58.5355}.

\bibitem[Katanaev(2011)]{adiab}
Katanaev, M.~O.
\newblock Adiabatic theorem for finite dimensional quantum mechanical systems.
\newblock \emph{Russian Physics Journal}, 54\penalty0 (3):\penalty0 342–353,
  August 2011.
\newblock ISSN 1573-9228.
\newblock \doi{10.1007/s11182-011-9620-5}.
\newblock URL \url{http://dx.doi.org/10.1007/s11182-011-9620-5}.

\bibitem[Katzgraber(2018)]{vanilla}
Katzgraber, H.~G.
\newblock Viewing vanilla quantum annealing through spin glasses.
\newblock \emph{Quantum Science and Technology}, 3\penalty0 (3):\penalty0
  030505, June 2018.
\newblock ISSN 2058-9565.
\newblock \doi{10.1088/2058-9565/aab6ba}.
\newblock URL \url{http://dx.doi.org/10.1088/2058-9565/aab6ba}.

\bibitem[Lucas(2014)]{qubo2}
Lucas, A.
\newblock Ising formulations of many np problems.
\newblock \emph{Frontiers in Physics}, 2, 2014.
\newblock ISSN 2296-424X.
\newblock \doi{10.3389/fphy.2014.00005}.
\newblock URL \url{http://dx.doi.org/10.3389/fphy.2014.00005}.

\bibitem[McGeoch \& Farré(2021)McGeoch and Farré]{adv}
McGeoch, C. and Farré, P.
\newblock The advantage system: Performance update, technical report.
\newblock Webpage, 2021.
\newblock URL
  \url{https://www.dwavesys.com/media/kjtlcemb/14-1054a-a_advantage_system_performance_update.pdf}.
\newblock Accessed on 2024-03-08.

\bibitem[Neukart et~al.(2017)Neukart, Compostella, Seidel, von Dollen, Yarkoni,
  and Parney]{traffic}
Neukart, F., Compostella, G., Seidel, C., von Dollen, D., Yarkoni, S., and
  Parney, B.
\newblock Traffic flow optimization using a quantum annealer, 2017.

\bibitem[Pakhomchik et~al.(2022)Pakhomchik, Yudin, Perelshtein, Alekseyenko,
  and Yarkoni]{scheduling}
Pakhomchik, A.~I., Yudin, S., Perelshtein, M.~R., Alekseyenko, A., and Yarkoni,
  S.
\newblock Solving workflow scheduling problems with qubo modeling, 2022.

\bibitem[Palackal et~al.(2023)Palackal, Poggel, Wulff, Ehm, Lorenz, and
  Mendl]{veh}
Palackal, L., Poggel, B., Wulff, M., Ehm, H., Lorenz, J.~M., and Mendl, C.~B.
\newblock Quantum-assisted solution paths for the capacitated vehicle routing
  problem, 2023.

\bibitem[Raymond et~al.(2023)Raymond, Stevanovic, Bernoudy, Boothby, McGeoch,
  Berkley, Farré, Pasvolsky, and King]{hybrid}
Raymond, J., Stevanovic, R., Bernoudy, W., Boothby, K., McGeoch, C.~C.,
  Berkley, A.~J., Farré, P., Pasvolsky, J., and King, A.~D.
\newblock Hybrid quantum annealing for larger-than-qpu lattice-structured
  problems.
\newblock \emph{ACM Transactions on Quantum Computing}, 4\penalty0
  (3):\penalty0 1–30, April 2023.
\newblock ISSN 2643-6817.
\newblock \doi{10.1145/3579368}.
\newblock URL \url{http://dx.doi.org/10.1145/3579368}.

\bibitem[Rieffel et~al.(2014)Rieffel, Venturelli, O’Gorman, Do, Prystay, and
  Smelyanskiy]{emb}
Rieffel, E.~G., Venturelli, D., O’Gorman, B., Do, M.~B., Prystay, E.~M., and
  Smelyanskiy, V.~N.
\newblock A case study in programming a quantum annealer for hard operational
  planning problems.
\newblock \emph{Quantum Information Processing}, 14\penalty0 (1):\penalty0
  1–36, December 2014.
\newblock ISSN 1573-1332.
\newblock \doi{10.1007/s11128-014-0892-x}.
\newblock URL \url{http://dx.doi.org/10.1007/s11128-014-0892-x}.

\bibitem[Stollenwerk et~al.(2018)Stollenwerk, Lobe, and Jung]{gate}
Stollenwerk, T., Lobe, E., and Jung, M.
\newblock Flight gate assignment with a quantum annealer, 2018.

\bibitem[Yan \& Sinitsyn(2022{\natexlab{a}})Yan and Sinitsyn]{nature_nonad}
Yan, B. and Sinitsyn, N.~A.
\newblock Analytical solution for nonadiabatic quantum annealing to arbitrary
  ising spin hamiltonian.
\newblock \emph{Nature Communications}, 13\penalty0 (1), April
  2022{\natexlab{a}}.
\newblock ISSN 2041-1723.
\newblock \doi{10.1038/s41467-022-29887-0}.
\newblock URL \url{http://dx.doi.org/10.1038/s41467-022-29887-0}.

\bibitem[Yan \& Sinitsyn(2022{\natexlab{b}})Yan and Sinitsyn]{qaoutperformance}
Yan, B. and Sinitsyn, N.~A.
\newblock Analytical solution for nonadiabatic quantum annealing to arbitrary
  ising spin hamiltonian.
\newblock \emph{Nature Communications}, 13\penalty0 (1), April
  2022{\natexlab{b}}.
\newblock ISSN 2041-1723.
\newblock \doi{10.1038/s41467-022-29887-0}.
\newblock URL \url{http://dx.doi.org/10.1038/s41467-022-29887-0}.

\end{thebibliography}
\bibliographystyle{icml2023}

\onecolumn\newpage
\section{Appendix}
We report here the tables with the parameter details of the synthetic problem settings used in the experiments.

\begin{table*}[!h]
\centering
\tiny
\setlength\tabcolsep{2pt}
\begin{tabularx}{\textwidth}{|l|l|l|l|X|X|}
\hline
\textbf{Instance} & \textbf{\# deliveries} & \textbf{battery budget} & \textbf{distribution} & \textbf{costs}                                       & \textbf{time intervals}                                                                                              \\ \hline
1                 & 8                      & 100                     & gaussian              & {[}55.2, 54.2, 41.1, 71.7, 55.2, 69.3, 46.0, 39.5{]} & {[}{[}15, 16{]}, {[}15, 18{]}, {[}10, 12{]}, {[}14, 15{]}, {[}10, 11{]}, {[}18, 20{]}, {[}11, 12{]}, {[}9, 12{]}{]}  \\ \hline
2                 & 8                      & 70                      & uniform               & {[}20.8, 5.5, 19.9, 26.2, 4.1, 2.7, 32.9, 61.8{]}    & {[}{[}19, 20{]}, {[}18, 20{]}, {[}16, 17{]}, {[}9, 10{]}, {[}17, 20{]}, {[}14, 15{]}, {[}16, 17{]}, {[}13, 16{]}{]}  \\ \hline
3                 & 8                      & 70                      & gaussian              & {[}38.6, 29.3, 25.2, 28.7, 27.2, 45.6, 47.3, 37.4{]} & {[}{[}8, 9{]}, {[}8, 11{]}, {[}11, 12{]}, {[}14, 15{]}, {[}10, 11{]}, {[}15, 18{]}, {[}8, 11{]}, {[}19, 20{]}{]}     \\ \hline
4                 & 7                      & 100                     & gaussian              & {[}57.2, 51.2, 45.6, 39.5, 44.2, 61.3, 44.0{]}       & {[}{[}19, 20{]}, {[}18, 20{]}, {[}13, 15{]}, {[}16, 19{]}, {[}9, 10{]}, {[}9, 11{]}, {[}14, 16{]}{]}                 \\ \hline
5                 & 8                      & 100                     & uniform               & {[}8.2, 86.8, 78.2, 61.0, 10.9, 49.9, 81.8, 44.5{]}  & {[}{[}19, 20{]}, {[}9, 10{]}, {[}11, 14{]}, {[}12, 13{]}, {[}18, 20{]}, {[}8, 11{]}, {[}10, 12{]}, {[}9, 12{]}{]}    \\ \hline
6                 & 8                      & 50                      & gaussian              & {[}25.4, 22.0, 50.0, 36.4, 40.3, 16.8, 19.3, 13.5{]} & {[}{[}16, 18{]}, {[}14, 16{]}, {[}12, 14{]}, {[}11, 14{]}, {[}14, 15{]}, {[}14, 15{]}, {[}11, 14{]}, {[}10, 11{]}{]} \\ \hline
7                 & 7                      & 50                      & gaussian              & {[}19.0, 29.2, 28.6, 29.9, 32.3, 19.3, 17.4{]}       & {[}{[}11, 13{]}, {[}12, 15{]}, {[}14, 15{]}, {[}15, 16{]}, {[}10, 13{]}, {[}17, 19{]}, {[}9, 10{]}{]}                \\ \hline
8                 & 8                      & 100                     & gaussian              & {[}47.0, 31.5, 54.1, 62.8, 54.9, 68.9, 54.0, 45.8{]} & {[}{[}18, 20{]}, {[}9, 11{]}, {[}15, 18{]}, {[}10, 13{]}, {[}18, 20{]}, {[}12, 15{]}, {[}16, 19{]}, {[}12, 14{]}{]}  \\ \hline
9                 & 8                      & 100                     & uniform               & {[}25.8, 31.5, 84.1, 94.5, 64.3, 47.2, 51.9, 96.6{]} & {[}{[}9, 11{]}, {[}13, 14{]}, {[}18, 20{]}, {[}12, 14{]}, {[}19, 20{]}, {[}15, 16{]}, {[}15, 18{]}, {[}8, 11{]}{]}   \\ \hline
10                & 7                      & 70                      & gaussian              & {[}42.8, 25.5, 37.9, 35.5, 18.8, 22.8, 34.7{]}       & {[}{[}10, 13{]}, {[}17, 18{]}, {[}13, 16{]}, {[}16, 19{]}, {[}11, 13{]}, {[}16, 17{]}, {[}17, 19{]}{]}               \\ \hline
11                & 8                      & 50                      & gaussian              & {[}11.4, 45.1, 12.8, 26.8, 26.7, 30.8, 28.1, 27.5{]} & {[}{[}14, 17{]}, {[}11, 13{]}, {[}19, 20{]}, {[}17, 20{]}, {[}14, 15{]}, {[}15, 17{]}, {[}11, 12{]}, {[}9, 12{]}{]}  \\ \hline
12                & 7                      & 100                     & gaussian              & {[}63.4, 49.4, 35.0, 62.8, 46.1, 49.3, 43.8{]}       & {[}{[}12, 14{]}, {[}11, 14{]}, {[}9, 12{]}, {[}8, 11{]}, {[}8, 11{]}, {[}19, 20{]}, {[}17, 19{]}{]}                  \\ \hline
\end{tabularx}

\caption{Sampled problem sets. Drones are fixed to be $10$, deliveries $N$ are randomly picked in $\{7,8\}$ and battery budgets in $\{50, 70, 100\}$. Also, a distribution is randomly chosen between uniform and Gaussian for each instance, and according to this, we generate $N$-dimensional lists of delivery energy costs and delivery time intervals, respectively. Time windows are generated from 8 a.m. to 8 p.m. and are at least one and at most two hours long, while costs $c_j$ are required to be smaller than the budget.}
\label{tab:N78_bench}
\end{table*}

\begin{table*}[!h]
\centering
\small
\setlength\tabcolsep{2pt}
\begin{tabularx}{\textwidth}{|l|l|l|l|X|X|}
\hline
\textbf{Instance} & \textbf{\# deliveries} & \textbf{battery budget} & \textbf{distribution} & \textbf{costs}                                       & \textbf{time intervals}                                                                                           \\ \hline
1                 & 4                      & 50                      & uniform               & {[}13.3, 41.9, 43.7, 19.0{]}                         & {[}{[}16, 18{]}, {[}13, 14{]}, {[}10, 11{]}, {[}13, 16{]}{]}                                                      \\ \hline
2                 & 5                      & 70                      & gaussian              & {[}50.8, 37.8, 32.5, 21.1, 25.6{]}                   & {[}{[}10, 11{]}, {[}14, 15{]}, {[}17, 19{]}, {[}9, 10{]}, {[}13, 16{]}{]}                                         \\ \hline
3                 & 6                      & 70                      & gaussian              & {[}52.2, 28.9, 41.9, 41.7, 36.2, 44.3{]}             & {[}{[}8, 9{]}, {[}18, 19{]}, {[}10, 11{]}, {[}12, 14{]}, {[}17, 19{]}, {[}12, 14{]}{]}                            \\ \hline
4                 & 7                      & 50                      & uniform               & {[}17.2, 40.1, 10.3, 30.4, 26.4, 40.5, 15.8{]}       & {[}{[}14, 17{]}, {[}16, 18{]}, {[}18, 20{]}, {[}16, 19{]}, {[}14, 16{]}, {[}16, 17{]}, {[}18, 19{]}{]}            \\ \hline
5                 & 8                      & 100                     & gaussian              & {[}38.9, 42.8, 52.8, 50.6, 34.1, 52.0, 66.5, 54.5{]} & {[}{[}15, 18{]}, {[}12, 13{]}, {[}13, 15{]}, {[}16, 18{]}, {[}9, 12{]}, {[}15, 17{]}, {[}17, 19{]}, {[}8, 9{]}{]} \\ \hline
\end{tabularx}
\caption{Sampled problem sets. Drones are fixed to be $10$, deliveries $N$ are considered to increase in $\{4,5,6,7,8\}$ and battery budgets in $\{50, 70, 100\}$. Also, a distribution is randomly chosen between uniform and Gaussian for each instance, and according to this, we generate $N$-dimensional lists of delivery energy costs and delivery time intervals, respectively. Time windows are generated from 8 a.m. to 8 p.m. and are at least one and at most two hours long, while costs $c_j$ are required to be smaller than the budget.}
\label{tab:incr_N_bench}
\end{table*}

\end{document}